%% file: ex_article.tex
\documentclass[hidelinks,onefignum,onetabnum]{siamart251216}

\input{ex_shared}

\ifpdf
\hypersetup{
  pdftitle={Symplectic ROM for Hamiltonian},
  pdfauthor={Y. Chen, et al.}
}
\fi

\begin{document}

\maketitle

\begin{abstract}
We introduce a novel data-driven symplectic reduced‐order modeling (ROM) framework for high‐dimensional Hamiltonian systems that unifies latent‐space discovery and dynamics learning within a single, end‐to‐end neural architecture. The encoder–decoder is built from H\'enon neural networks (H\'enonNets) and may be augmented with linear $G$-reflector layers, this yields an exact symplectic map between full and latent phase spaces. Latent dynamics are advanced by a symplectic flow map implemented as a H\'enonNet.
This unified neural architecture ensures exact preservation of the underlying symplectic structure at the reduced-order level, significantly enhancing the fidelity and long-term stability of the resulting ROM.  We validate our method through comprehensive numerical experiments on canonical Hamiltonian systems. The results demonstrate the method's capability for accurate trajectory reconstruction, robust predictive performance beyond the training horizon, and accurate Hamiltonian preservation. These promising outcomes underscore the effectiveness and potential applicability of our symplectic ROM framework for complex dynamical systems across a broad range of scientific and engineering disciplines. {The implementation associated with this work is publicly available at \url{https://github.com/chenyongssss/Symp_NN}.}
\end{abstract}




\section{Introduction}  \label{sec:introduction}

In this work, we investigate symplectic reduced-order modeling (ROM) for parametric Hamiltonian systems, following the framework introduced by Peng et al. \cite{peng2016}. Parametric Hamiltonian systems with high-dimensional state spaces are ubiquitous across various scientific and engineering domains, including molecular dynamics \cite{leimkuhler2004simulating}, celestial mechanics \cite{hairer2006geometric}, plasma physics \cite{morrison1998hamiltonian}, and quantum systems \cite{lubich2008quantum}. However, executing high-fidelity simulations for these systems can be computationally prohibitive. ROM techniques offer powerful alternatives for mitigating these computational costs by leveraging the underlying assumption that solutions, despite their apparent high dimensionality, frequently reside on low-dimensional manifolds embedded within the high-dimensional ambient space \cite{benner2015survey, quarteroni2015reduced}.

ROM techniques are generally classified into linear and nonlinear approaches. Linear ROM methods have been extensively developed and applied, particularly for linear systems. Modal analysis \cite{rouleau2017comparison, bonvin1982unified}, which uses eigenmode decomposition to represent dominant system dynamics, and proper orthogonal decomposition (POD) coupled with Galerkin projection \cite{moore1981principal, parrilo1999model,kunisch2002galerkin}, which identifies optimal linear subspaces for approximation, are foundational approaches. Krylov subspace methods \cite{odabasioglu2003prima, freund2003model, bai2002krylov} have also demonstrated their efficacy in reducing large-scale linear systems.

Nonetheless, real-world applications frequently exhibit pronounced nonlinear behavior, necessitating nonlinear ROM techniques. Recent advancements include the discrete empirical interpolation method (DEIM) \cite{chaturantabut2010nonlinear, drohmann2012reduced,drmac2016new}, which enhances POD-based approaches to handle nonlinearities efficiently. Manifold learning strategies \cite{peng2014online, millan2013nonlinear}, leveraging locally linear embeddings and kernel methods, have been developed to capture complex nonlinear features. Additionally, recent machine learning (ML) methodologies, notably deep autoencoders \cite{gonzalez2018deep, pagani2022statistical, pant2021deep} and physics-informed losses \cite{qian2020lift, chen2021physics}, have shown considerable promise in modeling nonlinear dynamics while integrating physical constraints.

Particularly relevant are our recent efforts in structure-preserving ML techniques designed for learning surrogates that can predict long-time dynamics. \cite{linot2023stabilized} introduced a stabilized neural ODE architecture that decomposes dynamics into sparse linear and dense nonlinear components, effectively capturing shocks and chaotic attractors. \cite{serino2025fast} proposed a fast–slow neural network that enforces a constrained invariant slow manifold, enabling accurate extrapolation of singularly perturbed systems. \cite{loya2025structure} developed structure-preserving NODEs with exponential integrators, offering provable Lyapunov stability for stiff systems. These efforts underscore the value of embedding physical and geometric structures into ML surrogates, motivating our current development of a symplectic ML-based ROM framework.


Recent years have also seen substantial progress in neural architectures that preserve Hamiltonian or symplectic structure at the level of learned dynamics. Representative examples include Hamiltonian neural networks \cite{greydanus2019hamiltonian}, which learn dynamics through an inferred Hamiltonian, and symplectic neural architectures such as SympNets \cite{jin2020sympnets}, together with subsequent structure-preserving symplectic map approximators \cite{duruisseaux2023approximation}. Closely related from the model-reduction perspective are recent latent-space approaches that impose variational or thermodynamic consistency, including variationally consistent Hamiltonian model reduction \cite{gruber2025variationally} and thermodynamically consistent latent dynamics identification for parametric systems \cite{he2025thermodynamically}. More broadly, recent work has explored reversible, irreversible, and metriplectic structure-preserving neural architectures beyond the ROM setting \cite{gruber2023reversible,gruber2025efficiently}. The present paper situates itself within this literature, focusing on the development of symplectic reduced-order models for parametric Hamiltonian systems.

Despite the recent advances, ROM techniques for parametric Hamiltonian systems face two primary limitations. Firstly, conventional ROM approaches often fail to preserve the inherent symplectic structure of Hamiltonian systems, potentially resulting in numerical instability. Secondly, due to the non-dissipative nature of Hamiltonian dynamics, global linear subspace methods frequently become ineffective, largely attributed to the slow decay of the Kolmogorov $n$-width associated with the solution manifold. Recent efforts addressing these challenges include the symplectic manifold Galerkin (SMG) projection method \cite{buchfink2023symplectic}, which generalizes structure-preserving model reduction from linear subspaces (see \cite{peng2016}) to symplectic nonlinear mappings. However, the initial implementation of SMG employed weakly symplectic mappings constructed through deep convolutional autoencoders, which lacked interpretability and required substantial hyperparameter tuning. A subsequent advancement \cite{sharma2023symplectic} introduced quadratic approximations within the SMG framework, effectively applied to nonlinear wave equations but constrained to quadratic nonlinearities and computationally demanding.

To address these limitations, we propose a data-driven symplectic ROM built from H\'enon neural networks (H\'enonNets) \cite{burby2020fast, duruisseaux2023approximation}, optionally augmented with linear \(G\)-reflector layers \cite{mackey2004g,mackey2003determinant}. H\'enonNets serve as the primary nonlinear building blocks because each layer is an exact symplectic map and compositions remain symplectic. In addition, H\'enonNet-based architectures possess a universal approximation property for symplectic maps, which provides the theoretical basis for their use in nonlinear symplectic embeddings and latent dynamics learning, and underlies the universal approximation result for the composite symplectic embedding architecture established in Section~\ref{sec:architecture}.
The optional $G$-reflector layers provide lightweight, parameter-efficient linear symplectic transformations that can capture the dominant linear structure and yield faster convergence~\cite{mackey2003determinant}. Our framework comprises two tightly integrated components: \emph{symplectic latent-space discovery} and \emph{structure-preserving latent dynamics learning}. First, a structure-preserving autoencoder maps the high-dimensional phase space onto a low-dimensional latent manifold using a H\'enonNet-based encoder–decoder, with optional $G$-reflector layers to provide linear symplectic corrections. Second, latent dynamics are advanced by a symplectic flow map realized as a H\'enonNet. Because symplecticity is preserved under composition, the overall mapping remains exactly symplectic, while the H\'enonNet backbone endows the model with strong approximation power for nonlinear Hamiltonian flows.

For training, we implement a unified training framework with a multi-objective loss function that simultaneously enforces state reconstruction accuracy, latent-space prediction fidelity, and Hamiltonian conservation. This joint optimization allows the latent space to be dynamically informed, enabling the encoder to extract features optimal for both reconstruction and temporal evolution prediction. Additionally, we incorporate training strategies including multi-step autoregressive prediction and small noise injection to enhance long-term stability and robustness.

In summary, our architecture preserves the symplectic structure by construction and promotes Hamiltonian consistency and long-time numerical stability through the proposed nonlinear transformations and training framework. The proposed approach is intended to improve the accuracy and structural fidelity of ROMs for Hamiltonian systems, particularly when nonlinear effects dominate the underlying dynamics.

 The remainder of this paper is organized as follows. Section~\ref{sec:governing} reviews the symplectic geometry and governing Hamiltonian equations. Section~\ref{sec:architecture} details the H\'enonNet, nonlinear symplectic lifts, G-reflector constructions, and the composite embedding architecture with universal approximation guarantees. Section~\ref{sec:methods} formulates the data-driven symplectic ROM, including latent space discovery, latent dynamics learning, and the unified training procedure. Numerical experiments and performance comparisons appear in Section~\ref{sec:results}, and concluding remarks and future directions are provided in Section~\ref{sec:conclusion}.

\section{Dynamical systems and symplectic embedding} \label{sec:governing}


The seminal work by Peng et al.\ \cite{peng2016} introduced key definitions related to Hamiltonian systems and its compression. Before presenting our neural network–based enhancement, we first review several important concepts from \cite{peng2016} to lay the foundation for our approach.

Let $(\mathbb{V}, \Omega)$ denote a \emph{symplectic vector space}, where $\mathbb{V}$ is a real vector space of dimension $2n$, and $\Omega(\cdot, \cdot)$ is a nondegenerate, alternating bilinear form, termed the \emph{symplectic form}, defined by:
\begin{equation}
\Omega: \mathbb{V} \times \mathbb{V} \to \mathbb{R}.
\end{equation}
A critical property of symplectic spaces is the existence of a canonical basis $\{e_1, \dots, e_n \\,  w_1, \dots, w_n\}$ satisfying:

\begin{equation}
\Omega(e_i,e_j) \;=\; \Omega(w_i,w_j) \;=\; 0, 
\quad 
\Omega(e_i,w_j) \;=\; \delta_{ij}, 
\quad i,j = 1,\dots,n,
\label{eq:canonical-basis}
\end{equation}
where \(\delta_{ij}\) is the Kronecker delta function.

In this paper, we assume $\mathbb{V}$ is defined over the field $\mathbb{R}$, thereby identifying $\mathbb{V}$ with $\mathbb{R}^{2n}$. Under canonical coordinates $(q_1,\dots,q_n,p_1,\dots,p_n)$, the symplectic form $\Omega$ can be expressed canonically as:
\begin{equation}
\Omega = \sum_{i=1}^n dq_i \wedge dp_i.
\label{eq:canonical-symp-form}
\end{equation}
Moreover, the symplectic form $\Omega$ can be represented via the standard Poisson matrix $J_{2n}$ as:
\begin{equation}
\Omega(\mathbf{v}_1,\mathbf{v}_2) \;=\; \mathbf{v}_1^\top\,J_{2n}\,\mathbf{v}_2, 
\quad \forall\,\mathbf{v}_1,\mathbf{v}_2 \,\in\, \mathbb{V},
\label{eq:Omega-matrix}
\end{equation}
where
\begin{equation}
J_{2n} \;=\; 
\begin{bmatrix}
\,0_n & I_n\,\\
\,-I_n & 0_n\,
\end{bmatrix}.
\end{equation}
Here, $I_n$ and $0_n$ represent the identity and zero matrices of size $n \times n$, respectively. The matrix $J_{2n}$ satisfies
\begin{equation}
J_{2n}\,J_{2n}^\top \;=\; J_{2n}^\top\,J_{2n} \;=\; I_{2n},
\quad
J_{2n}^2 \;=\; -I_{2n}.
\end{equation}

We define the following two fundamental symplectic matrix manifolds used extensively throughout the paper:
\begin{equation}
\mathrm{Sp}(2n,\mathbb{R}) := \bigl\{S\in\mathbb{R}^{2n\times2n}\mid S^{\mathsf T} J_{2n}S=J_{2n}\bigr\},
\end{equation}
and
\begin{equation}
\mathrm{SpSt}(2k,2n) := \bigl\{U\in\mathbb{R}^{2n\times2k}\mid U^{\mathsf T}J_{2n}U=J_{2k}\bigr\},
\end{equation}
where $k\leq n$. The former is known as the \emph{real symplectic group}, consisting of all linear transformations preserving the symplectic structure. The latter, $\mathrm{SpSt}(2k,2n)$, is the \emph{symplectic Stiefel manifold}, which represents the set of all full-rank matrices whose columns form a symplectic basis for a $2k$-dimensional subspace within a larger $2n$-dimensional space. In the full-dimensional case \(k=n\), this definition reduces to
\begin{equation}
\mathrm{SpSt}(2n,2n)=\mathrm{Sp}(2n,\mathbb{R}),
\end{equation}
so the symplectic Stiefel manifold may be viewed as a natural extension of the real symplectic group to rectangular symplectic embeddings.

Consider a  Hamiltonian function $H: \mathbb{V} \to \mathbb{R}$. The dynamics of an autonomous Hamiltonian system on $\mathbb{V}$ are governed by \emph{Hamiltonian equations}:
\begin{equation}
\dot{\mathbf{q}} = \nabla_{\mathbf{p}}H(\mathbf{q},\mathbf{p}), \quad \dot{\mathbf{p}} = -\nabla_{\mathbf{q}}H(\mathbf{q},\mathbf{p}),
\label{eq:hamilton-equations}
\end{equation}
where \(\mathbf{q}\in\mathbb{R}^n\) and \(\mathbf{p}\in\mathbb{R}^n\) 
represent the position and momentum vectors, respectively. Defining the phase space variable $\mathbf{x} = [\mathbf{q}^\top \, \mathbf{p}^\top]^\top$, Hamiltonian equations can be succinctly expressed as:
\begin{equation}
\dot{\mathbf{x}} = J_{2n}\nabla_{\mathbf{x}}H(\mathbf{x}).
\label{eq:hamilton-compact}
\end{equation}

A fundamental property of Hamiltonian systems is the generation of \emph{symplectomorphisms} by its flow, preserving the symplectic form $\Omega$, the Hamiltonian $H$, and the volume of phase space. Preserving these invariants is crucial for maintaining the qualitative behaviors of Hamiltonian systems, which motivates the development and application of symplectic numerical integrators such as the St{\"o}rmer–Verlet method \cite{hairer2006geometric}. The method is given explicitly as:
\begin{equation}
\begin{aligned}
\mathbf{p}_{n+\tfrac{1}{2}} 
&=\; 
\mathbf{p}_n \;-\; 
\tfrac{\Delta t}{2}\,\nabla_{\mathbf{q}}\,H\bigl(\mathbf{q}_n,\mathbf{p}_{n+\tfrac{1}{2}}\bigr),
\\
\mathbf{q}_{n+1} 
&=\;
\mathbf{q}_n \;+\; 
\tfrac{\Delta t}{2}\,\Bigl(\nabla_{\mathbf{p}}\,H\bigl(\mathbf{q}_n,\mathbf{p}_{n+\tfrac{1}{2}}\bigr)
\;+\;
\nabla_{\mathbf{p}}\,H\bigl(\mathbf{q}_{n+1},\mathbf{p}_{n+\tfrac{1}{2}}\bigr)\Bigr),
\\
\mathbf{p}_{n+1} 
&=\;
\mathbf{p}_{n+\tfrac{1}{2}} \;-\; 
\tfrac{\Delta t}{2}\,\nabla_{\mathbf{q}}\,H\bigl(\mathbf{q}_{n+1},\mathbf{p}_{n+\tfrac{1}{2}}\bigr).
\end{aligned}
\label{eq:stormer-verlet}
\end{equation}

For the model reduction task of Hamiltonian systems, the goal is to identify a lower-dimensional latent space to approximate the high-dimensional system dynamics effectively. To maintain stability and physical consistency, this reduction must preserve the symplectic structure. Formally, consider two symplectic vector spaces $(\mathbb{V},\Omega)$ and $(\mathbb{W},\omega)$ with dimensions $2n$ and $2k$, respectively, with $k \ll n$. Each space has its symplectic form represented by: 
\begin{equation}
\Omega(\mathbf{v}_1,\mathbf{v}_2)=\mathbf{v}_1^\top J_{2n}\mathbf{v}_2,\quad \omega(\mathbf{z}_1,\mathbf{z}_2)=\mathbf{z}_1^\top J_{2k}\mathbf{z}_2,
\end{equation}
for vectors $\mathbf{v}_1,\mathbf{v}_2\in\mathbb{V}$ and $\mathbf{z}_1,\mathbf{z}_2\in\mathbb{W}$.

A smooth embedding $g: \mathbb{W}\to\mathbb{V}$ that relates the latent variables $\mathbf{z}\in\mathbb{W}$ to the full phase-space variables $\mathbf{v}\in\mathbb{V}$ is called a \emph{symplectic embedding (lift)} if it exactly preserves the symplectic form. Specifically, this means that the symplectic form on the latent space is the \emph{pullback} of the symplectic form on the full space via $g$:
\begin{equation}
g^*\Omega = \omega.
\end{equation}

This condition explicitly implies that the Jacobian $Dg(\mathbf{z})\in\mathbb{R}^{2n\times 2k}$ of the embedding  must satisfy the symplectic preservation condition:
\begin{equation}
Dg(\mathbf{z})^\top J_{2n} \, Dg(\mathbf{z}) = J_{2k}, \quad \forall \mathbf{z}\in\mathbb{W}.
\end{equation}
This constraint ensures exact preservation of the geometric structure, critical for maintaining long-term stability and physical fidelity in ROM of Hamiltonian systems.

For the linear case,  which serves as the foundation for many classical approaches, the symplectic embedding takes the form:
\begin{equation}
g(\mathbf{z}) = A\mathbf{z}, \quad \mathbf{z}\in\mathbb{W},
\end{equation}
with $A\in\mathbb{R}^{2n\times 2k} \in \mathrm{SpSt}(2k,2n)$ satisfying the symplectic condition:
\begin{equation}
A^\top J_{2n} A = J_{2k}.
\label{eq:symplectic-lift-linear}
\end{equation}
The symplectic inverse of such a matrix $A$ is given by:
\begin{equation}
A^+ = J_{2k}^\top A^\top J_{2n}.
\label{eq:symplectic-inverse}
\end{equation}
While the previous work \cite{peng2016} focuses on linear mappings, we will show that our approach generalizes these techniques by introducing nonlinear mappings to more effectively capture the intricate behaviors characteristic of realistic Hamiltonian systems.

\section{Symplectic neural network and embedding}
\label{sec:architecture}

In this section, we introduce a unified neural network architecture capable of approximating arbitrary symplectic maps primarily through compositions of nonlinear symplectic lifts constructed from H\'enonNets, with additional linear corrections provided by symplectic transformations constructed from G-reflectors. This framework, centered on advanced nonlinear mappings, lays the foundation for the model reduction method presented in Section~\ref{sec:methods}.


\subsection{H\'enonNets: symplectic neural networks}
In this subsection, we introduce the H\'enonNet architecture, a cornerstone of our proposed symplectic ROM framework. Initially developed in \cite{burby2020fast}, 
H\'enonNet constructs symplectic neural networks through structured compositions of elementary nonlinear symplectic mappings, capable of capturing complex Hamiltonian dynamics.

The fundamental building block is the H\'enon mapping, which enables intricate nonlinear transformations while preserving symplectic structure:
\begin{definition}[H\'enon Mapping]\label{def:Henon-map}
For a smooth scalar-valued function $V: \mathbb{R}^n \to \mathbb{R}$ and a constant vector $\mathbf{\eta} \in \mathbb{R}^n$, the H\'enon mapping $H(V,\eta): \mathbb{R}^{2n} \to \mathbb{R}^{2n}$ is defined as:
\begin{equation}
H(V,\eta)\begin{pmatrix}
\xv \\ \yv
\end{pmatrix} = \begin{pmatrix}
\yv+ \mathbf{\eta} \\ \xv + \nabla V(\yv)
\end{pmatrix},
\end{equation}
where $\nabla V$ represents the gradient of $V$.
\end{definition}

By utilizing neural networks to parameterize $V$, we construct the H\'enon layer, 
 as the fourth-power composition of a single H\'enon mapping:
\begin{equation}\label{eq:Henon-layer}
\mathcal{H}(V,\eta)
=
\underbrace{H(V,\eta)\circ H(V,\eta)\circ H(V,\eta)\circ H(V,\eta)}_{\text{fourth-power composition}}.
\end{equation}
For brevity, we may also write \(\mathcal{H}(V,\eta)=[H(V,\eta)]^{4}\), where the superscript \(4\) denotes repeated composition rather than an elementwise fourth power.
The complete H\'enonNet architecture is built through a sequential composition of multiple H\'enon layers:
\begin{definition}[H\'enonNet]\label{def:HenonNet}
A H\'enonNet with depth $N$ consists of:
\begin{itemize}
    \item A sequence of fully connected neural networks $\{V_i\}_{i=1}^N$, where each $V_i: \mathbb{R}^n \to \mathbb{R}$
    \item A sequence of trainable vectors $\{\eta_i\}_{i=1}^N \subset \mathbb{R}^n$
\end{itemize}
The H\'enonNet transformation is constructed through the sequential composition:
\begin{equation}
    \text{H\'enonNet} = \mathcal{H}(V_N,\eta_N) \circ \cdots \circ \mathcal{H}(V_1,\eta_1)
\end{equation}
\end{definition}

A crucial property of H\'enonNets is their explicit invertibility, which we formally establish below and will be leveraged in this work:
\begin{proposition}
 Every H\'enonNet transformation admits an analytical inverse. Specifically, for the H\'enon mapping $H(V,\eta)$, the inverse mapping $H^{-1}(V,\eta): \mathbb{R}^{2n} \to \mathbb{R}^{2n}$ is given by:
\begin{equation}
    H^{-1}(V,\eta)\begin{pmatrix}
\mathbf{u} \\ \mathbf{v}
\end{pmatrix} = \begin{pmatrix}
\mathbf{v} - \nabla V(\mathbf{u} - \mathbf{\eta}) \\ \mathbf{u} - \mathbf{\eta}
\end{pmatrix}.\end{equation}
\end{proposition}

Consequently, the complete H\'enonNet inverse is obtained through:
\begin{equation}
\text{H\'enonNet}^{-1} = \mathcal{H}^{-1}(V_1,\eta_1) \circ \cdots \circ \mathcal{H}^{-1}(V_N,\eta_N),\end{equation}
where each $\mathcal{H}^{-1}(V_i,\eta_i) = \left[H^{-1}(V_i,\eta_i)\right]^4$.
This explicit invertibility establishes H\'enon\-Nets as a special class of invertible neural networks~\cite{kobyzev2020normalizing}, with the non-trivial advantage that the inverse can be computed analytically without iterative procedures or additional computational overhead.
In fact, H\'enonNets can be considered as a symplectic generalization of the triangular normalizing flows~\cite{kobyzev2020normalizing}.


The fundamental property that makes H\'enonNet suitable for nonlinear Hamiltonian systems is its intrinsic preservation of symplectic structure, which is guaranteed independently of the neural network weights chosen for $V_i$. More critically, H\'enonNet has been proven to possess a universal approximation theorem for nonlinear symplectic mappings~\cite{turaev2002polynomial}, making it suitable for our symplectic model reduction framework. This universal approximation property enables H\'enonNets to capture nonlinear phase space geometries that are beyond the reach of traditional linear methods.



\subsection{Nonlinear symplectic embedding}
\label{sec:nonlinear}
Having established H\'enonNets for symplectic transformations on $\mathbb{R}^{2n}$, we now extend this framework to construct nonlinear symplectic embeddings that fundamentally transcend linear embedding limitations.

\begin{definition} A nonlinear {\bf symplectic embedding} from a vector space $\mathbb{W}\in \mathbb{R}^{2k}$ to a vector space $\mathbb{V} \in \mathbb{R}^{2n}$ can be defined as
\begin{align}\label{eq:nonlinear-sym-lift}
\begin{cases}
\hat{\xv} & = P \yv + \mathbf{\eta},\\
\hat{\yv} & = - P \xv + \grad V(P \yv).
\end{cases}
\end{align}
where $(\xv^\top,\yv^\top)^\top\in \mathbb{W}, (\hat{\xv}^\top,\hat{\yv}^\top)^\top\in \mathbb{V}$, $ P$ is a linear embedding map from the vector space of dim $k$ to the vector space of dim $n$, $V$ is an arbitrary smooth function in $\mathbb{V}$, and $\eta$ is a bias vector.
\end{definition}

We denote this nonlinear symplectic embedding as $g: \mathbb{W} \to \mathbb{V}$. The Jacobian of the map $g$ is:
\begin{align}
Dg : = \begin{bmatrix}
 0 & P \\
-P & H_V P
\end{bmatrix},
\end{align}
where $H_V$ standards for the Hessian matrix of the function $V$. 
 The following theorem characterizes the condition for symplectic preservation:
\begin{theorem}
The mapping $g$ preserves the symplectic form if and only if $P \in {\rm St}(k, n)$, i.e, the Stiefel manifold, explicitly:
\begin{equation}
(Dg)^\top J_{2n} Dg = J_{2k}.
\end{equation}
\end{theorem}

We naturally define the inverse mapping of the symplectic embedding, where we leverage the invertibility of H\'enonNets:
\begin{definition}The one-sided symplectic pseudo-inverse of the {\bf symplectic e\-mbedding} is defined by
\end{definition}
Note that the Jacobian of $g^{-1}$ is 
\begin{align}
D g^{-1} : = \begin{bmatrix}
 P^\top H_V & -P^\top \\
P^\top & 0
\end{bmatrix}
\end{align}
Note that the Jacobian satisfies the symplectic inverse definition
\begin{equation}
D g^{-1} = J_{2k}^\top\left(D g\right)^\top J_{2n}.
\end{equation}
This construction provides a natural nonlinear generalization of linear symplectic embeddings while maintaining the geometric properties required for Hamiltonian systems.

The universal approximation capability of H\'enonNets implies that any symplectic diffeomorphism, including the nonlinear embedding defined above, can be approximated arbitrarily closely by composing a fixed linear embedding with H\'enonNet transformations. Hence, in practical computations, it suffices to choose the simplest possible embedding map, the canonical inclusion map:
\begin{equation}\label{eq:insert}
    \iota = 
\begin{pmatrix}
I_{n\times k} & 0\\[3pt]
0 & I_{n\times k}
\end{pmatrix} \quad\in\mathrm{SpSt}(2k,2n),
\end{equation}
thus significantly simplifying implementation without sacrificing approximation accuracy.

Specifically, we define our nonlinear symplectic embedding in practice as the direct composition:
\begin{equation}\label{eq:direct-henon-embed}
  \text{H\'enonNet}\circ\iota  ((\xv^\top,\yv^\top)^\top))
  = \mathcal{H}(V_N,\eta_N)\circ\cdots\circ\mathcal{H}(V_1,\eta_1)\circ\iota ((\xv^\top,\yv^\top)^\top))
\end{equation}
where each \(\mathcal{H}(V_i,\eta_i)\) denotes a single H\'enon layer, which is defined previously in \eqref{eq:Henon-layer}, and \(\iota\) is fixed and not learned.

This simplified composition achieves two critical goals:
\begin{itemize}
  \item It leverages the universal nonlinear approximation power of H\'enonNets, capable of capturing intricate nonlinear structures inherent in Hamiltonian systems.
  \item It ensures strict preservation of the symplectic structure, essential for long-term stability and accuracy of the reduced-order model.
\end{itemize}

Thus, while general nonlinear symplectic embeddings may utilize learned linear embeddings \(P\), the practical equivalence due to universal approximation allows us to adopt a simpler, fixed embedding strategy. Throughout this paper, unless otherwise stated, the simplified composition (\ref{eq:direct-henon-embed}) is employed due to its practical efficiency and robust approximation capacity.

\subsection{Linear symplectic embedding}
\label{sec:linear}

Although it is not the main focus of the current study, in this subsection, we propose an alternative approach of the cotangent lift from the original symplectic ROM work \cite{peng2016}.
One issue associated with cotangent lift is its incomplete parameterization of $\mathrm{SpSt}(2k,2n)$, which shall be overcome by the complete parameterization through a linear symplectic layer proposed herein.
To this end, we introduce \emph{G-reflector} \cite{mackey2004g} layers, each implementing an exact, parameter-efficient linear symplectic map.


A G-reflector is defined for \(\mathbf{u}\in\mathbb{R}^{2n}\) and \(\beta\in\mathbb{R}\) by
\(
G \;=\; I_{2n} \;+\;\beta\,\mathbf{u}\,\mathbf{u}^{T}J_{2n},
\)
with its analytical inverse being
\(
G^{-1} = I_{2n} - \beta\,\mathbf{u}\,\mathbf{u}^{T}J_{2n}.
\)
It is easy to see \(G,\,G^{-1}\in\mathrm{Sp}(2n,\mathbb{R})\). 
The theoretical foundation rests on the result that any symplectic matrix can be factored into at most $4n$ such reflectors~\cite{mackey2003determinant}, providing a complete parameterization of the linear symplectic group. Each G-reflector introduces $2n+1$ trainable parameters, reducible to $2n$ under normalization of $\mathbf{u}$, yielding a computationally tractable representation. This layer can be easily incorporated into a ML framework with trainable $\mathbf{u}$ and $\beta$.


In a recent application to plasma physics \cite{Drimalas2025SymplecticNN}, we employed the G-reflector as the base layer to approximate a $4\times 4$ symplectic matrix.  The experiments confirmed the practical robustness of reflector layers in low-dimensional particle dynamics.
However, its full potential in a higher-dimension phase space has not been explored. 


To lift from latent spaces $\mathbb{R}^{2k}$ to ambient spaces $\mathbb{R}^{2n}$, we again use the rectangular inclusion \(\iota\) defined in \eqref{eq:insert}
and parameterize the full symplectic Stiefel embedding \(U\in\mathrm{SpSt}(2k,2n)\) by
\begin{equation}
U = S\,\iota,\
\end{equation}
where $S\;\in\;\mathrm{Sp}(2n,\mathbb{R})$ is realized through a small stack of G-reflectors.
For a formally complete parametrization, $4n$ G-reflector layers are needed. However, in practice, we find that adding more linear layers offers limited improvement in learning the overall dynamics, as linear compression is suboptimal for Hamiltonian systems and is rarely the primary bottleneck.


It should be noted that the mapping from $\mathrm{Sp}(2n,\mathbb{R})$ to the symplectic Stiefel manifold $\mathrm{SpSt}(2k,2n)$ through the construction $U = S\,\iota$ is a surjective  map \cite{bendokat2021real}. Consequently, the proposed parameterization of $S$ provides a complete way to represent every element of $\mathrm{SpSt}(2k,2n)$. In fact, a stronger result was discussed in \cite{bendokat2021real}: the manifold $\mathrm{SpSt}(2k, 2n)$ can be identified with the quotient space:
\begin{equation}
\mathrm{SpSt}(2k, 2n) \simeq \mathrm{Sp}(2n, \mathbb{R})/\mathrm{Sp}(2(n-k), \mathbb{R}).\end{equation}
where the quotient structure describes the redundancy in the parameterization arising from elements in $\mathrm{Sp}(2(n-k), \mathbb{R})$. 



The parameterization through $S$ enables efficient learning of the embedding while preserving geometric constraints throughout the optimization process. 
This can be also used as an optional linear layer integrated alongside the H\'enonNet layers. 
These linear corrections may capture primary linear dynamics for weakly linear cases, yielding faster convergence, and reduced parameter redundancy.


\subsection{Composite architecture and symplectic universal approximation}

We now assemble the previously introduced components into a unified composite architecture that achieves universal approximation capabilities for symplectic embeddings while maintaining structural guarantees. This composite framework integrates the nonlinear expressiveness of H\'enonNets with auxiliary linear symplectic corrections.


\begin{definition}[Composite Symplectic Embedding]\label{def:composite-embedding}
Given a symplectic embedding matrix $\iota \in \mathrm{SpSt}(2k, 2n)$, a H\'enonNet $\text{H}:\mathbb{R}^{2n}\to \mathbb{R}^{2n}$, and a linear symplectic transformation $G \in \mathrm{Sp}(2n, \mathbb{R})$ composed of $G$-reflectors, define the composite symplectic embedding as:
\begin{equation} \label{eq:composite}
    \sigma = \text{H} \circ G \circ \iota : \mathbb{R}^{2k} \to \mathbb{R}^{2n}.
\end{equation}

\end{definition}

The composite architecture admits a natural inverse construction that preserves symplectic structure:

\begin{definition}[Inverse Composite Embedding]\label{def:inverse-composite}
The inverse of the composite symplectic embedding is given by:
\begin{equation}
\sigma^{-1} = \iota^+ \circ G^{-1} \circ \text{H}^{-1} : \mathbb{R}^{2n} \to \mathbb{R}^{2k}
\end{equation}
where $\iota^+$ denotes the symplectic pseudo-inverse of $\iota$, $G^{-1}$ is the inverse of the $G$-reflector composition, and $\text{H}^{-1}$ is the inverse H\'enonNet.
\end{definition}

The symplectic pseudo-inverse $\iota^+$ is uniquely defined by the condition $\iota^+ \iota = I_{2k}$ and satisfies $(\iota^+)^\top J_{2n} \iota^+ = J_{2k}$, ensuring that the inverse mapping preserves symplectic structure on the embedded submanifold.

This composite structure possesses universal approximation capability for symplectic embeddings, as formalized in the following theorem:

\begin{theorem}[Universal Approximation for Symplectic Embeddings]\label{thm:universal-approx}
Let $\Phi:\mathbb{R}^{2k}\to \mathbb{R}^{2n}$ be an arbitrary symplectic embedding defined on a compact domain $\Omega\subset \mathbb{R}^{2k}$. For any $\varepsilon>0$, there exist a composition \(\mathcal{G}\) of \(G\)-reflectors and a H\'enonNet \(\mathrm{H}\) such that the composite map \(\sigma\) defined in \eqref{eq:composite} satisfies
\begin{equation}
\sup_{z\in\Omega} \|\Phi(z) - \sigma(z)\| < \varepsilon
\end{equation}
and
\begin{equation}
 \mathrm{D}\sigma(z)^\top J_{2n} \mathrm{D}\sigma(z) = J_{2k},
\end{equation}
where $\mathrm{D}\sigma(z)$ denotes the Jacobian matrix of $\sigma$ at $z$.
\end{theorem}

\begin{proof}
Since \(\Phi\) is a symplectic embedding, its image \(\Phi(\Omega)\) is a symplectic submanifold of \((\mathbb{R}^{2n},J_{2n})\). By a standard symplectic neighborhood theorem \cite[Theorem 3.22]{mcduff2017introduction}, there exists a local symplectomorphism \(\Psi\) defined near \(\iota(\Omega)\) such that
\begin{equation}
\Phi(z)=\Psi(\iota(z)), \qquad z\in\Omega.    
\end{equation}
Therefore, it suffices to approximate the ambient symplectic map \(\Psi\) on the compact set \(\iota(\Omega)\).

By the universal approximation property of H\'enonNets for symplectic maps \cite{duruisseaux2023approximation}, for any \(\varepsilon>0\), one can choose a H\'enonNet \(\mathrm{H}\) such that
\begin{equation}
\sup_{x\in \iota(\Omega)} \|\Psi(x)-\mathrm{H}(x)\|<\varepsilon.
\end{equation}
The optional linear correction \(\mathcal{G}\) may be taken as the identity, or more generally as a composition of \(G\)-reflectors when additional linear symplectic adjustment is desired. Hence, with
\begin{equation}
\sigma=\mathrm{H}\circ \mathcal{G}\circ \iota,
\end{equation}
we obtain the required uniform approximation on \(\Omega\).

Finally, \(\iota\) is symplectic on its image, \(\mathcal{G}\) is symplectic as a composition of \(G\)-reflectors, and \(\mathrm{H}\) is symplectic by construction. Since symplecticity is preserved under composition,
\begin{equation}
\mathrm{D}\sigma(z)^\top J_{2n}\mathrm{D}\sigma(z)=J_{2k},
\qquad z\in\Omega,
\end{equation}
which completes the argument.
\end{proof}

\begin{remark}
The depth of the H\'enonNet and the number of $G$-reflectors can be determined based on the desired approximation accuracy $\varepsilon$ and the complexity of the target symplectic mapping $\Phi$. The H\'enonNet depth primarily determines the nonlinear approximation capacity, while G-reflectors enhance numerical stability.
\end{remark}

This universal approximation result establishes that our composite architecture can approximate symplectic embeddings while maintaining exact symplectic structure preservation. While the embedding map \(\sigma\) itself is not globally invertible (as its domain and codomain have different dimensions), its inverse map \(\sigma^{-1}\) is well-defined and smooth on the embedded symplectic submanifold. Consequently, the composite architecture permits bidirectional symplectic transformations between the latent and ambient spaces.

Furthermore, the composite architecture allows for direct optimization of the embedding parameters using gradient-based methods. The combination of theoretical universality, structural preservation, and practical implementability makes this composite framework suitable for the symplectic autoencoder and ROM methodology developed in the next section.

\section{Symplectic model reduction}
\label{sec:methods}



\begin{figure}[htp]
    \centering
    \includegraphics[width=0.8\linewidth]{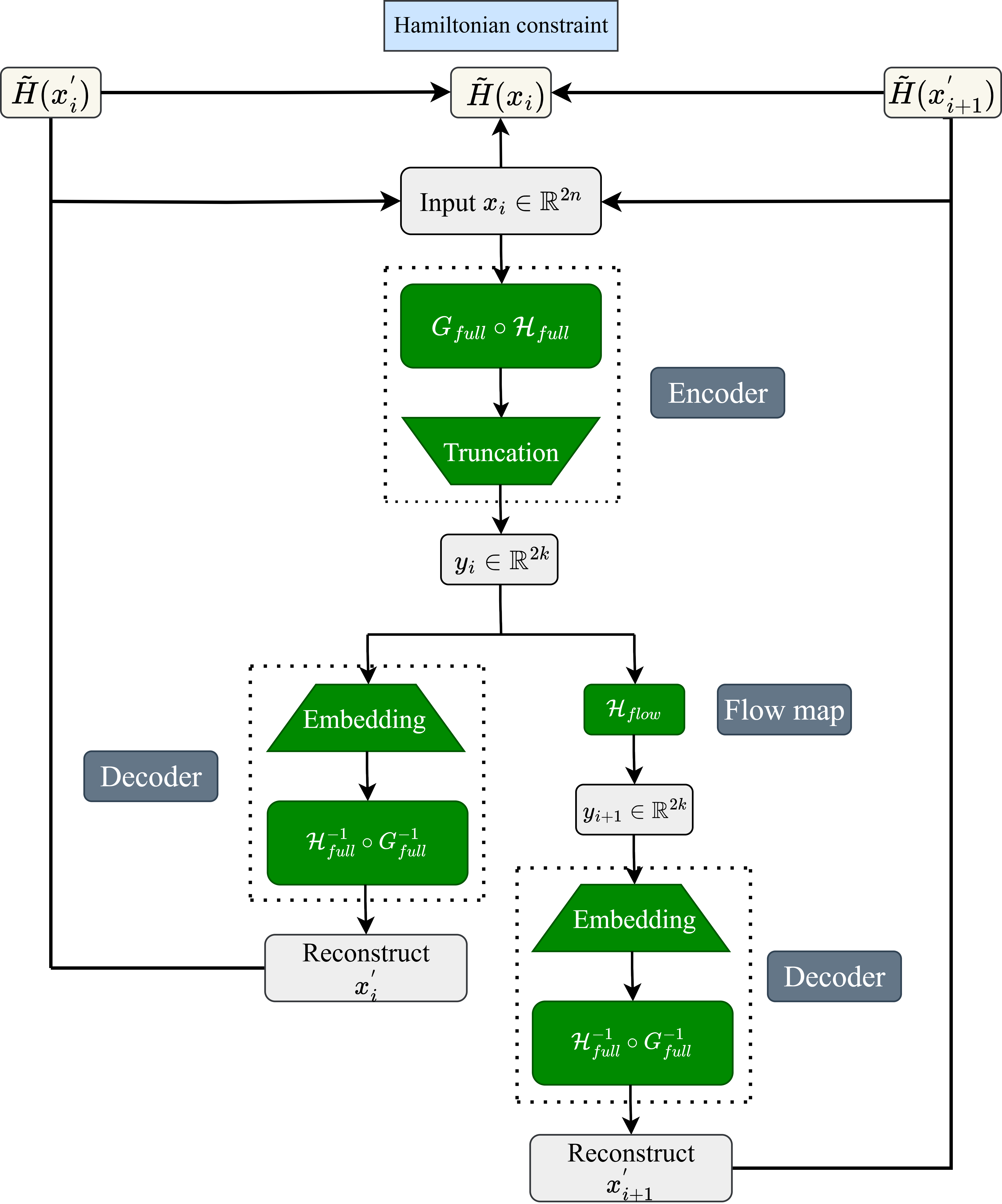}
    \caption{Schematic illustration of the symplectic neural network architecture for ROM of Hamiltonian systems. We note that the G-reflector in the autoencoder is only optional.}
    \label{fig:framework}
\end{figure}

In this section, we introduce a data-driven symplectic ROM framework for Hamiltonian systems. 
The framework consists of two main stages: (1) symplectic latent-space discovery, accomplished through a structure-preseving autoencoder that maps the high-dimensional phase space onto a lower-dimensional latent manifold; and (2) structure-preserving latent dynamics learning, employing a symplectic flow network to capture latent temporal evolution. These components are unified in a single neural architecture, as illustrated in Figure~\ref{fig:framework}.

\subsection{Latent space discovery}
\label{sec:composition}


The nonlinear encoder $f_{\mathrm{enc}}:\mathbb{R}^{2n}\to\mathbb{R}^{2k}$ is primarily constructed from a H\'enonNet (\(\mathcal{H}_{\mathrm{full}}\)). Optionally, to enhance stability and accelerate convergence—especially in cases where prominent linear structures exist—this nonlinear mapping can be supplemented by linear symplectic corrections composed of G-reflectors (\(G_{\mathrm{full}}\)). Formally, the encoder takes the form:
\begin{equation}
    \fenc(\mathbf{x}) = \tau\circ G_{\mathrm{full}}\circ \mathcal{H}_{\mathrm{full}}(\mathbf{x}), \mathbf{x}\in \mathbb{R}^{2n}
\end{equation}
The truncation operator $\tau: \mathbb{R}^{2n} \to \mathbb{R}^{2k}$ is defined as:
\begin{equation}
    \tau(\mathbf{z}) = [z_{1:k},\, z_{n+1:n+k}], \quad \mathbf{z}\in\mathbb{R}^{2n}.
\end{equation}
This choice ensures that $\tau$ acts as a symplectic subplace projection, and coincides with the symplectic pseudo-inverse of the embedding map $\iota$ defined previously in \eqref{eq:insert}. Consequently, the dimensionality reduction performed by the encoder rigorously preserves the symplectic structure.

The decoder $f_{\mathrm{dec}}:\mathbb{R}^{2k}\to\mathbb{R}^{2n}$ serves as the inverse operation, combining embedding matrix $\iota:\mathbb{R}^{2k}\to\mathbb{R}^{2n}$ with inverse mappings of the encoder components:
\begin{equation}
    \fdec(\mathbf{y}) = \mathcal{H}_{\mathrm{full}}^{-1}\circ G_{\mathrm{full}}^{-1}\circ \iota(\mathbf{y}),\mathbf{y}\in \mathbb{R}^{2k}
\end{equation}
where the embedding operator $\iota$ places $\mathbf{y}$ into the appropriate positions of a zero-initialized vector in $\mathbb{R}^{2n}$.

Given temporal snapshots of system states $\{\mathbf{x}_i\}_{i=0}^{T}\subset\mathbb{R}^{2n}$, we can train the autoencoder by minimizing the mean squared reconstruction error:
\begin{equation}
    \mathcal{L}_r = \frac{1}{T+1}\sum_{i=0}^{T}\|\xv_i - \fdec(\fenc(\xv_i))\|_2^2.
\end{equation}

\begin{remark}
The constructed encoder and decoder inherently preserve symplectic structure due to the composition of symplectic components, rendering the approach particularly effective for Hamiltonian systems.
\end{remark}

\subsection{Latent dynamics learning}


For the system states $\{\mathbf{x}_i\}_{i=0}^{T}\subset\mathbb{R}^{2n}$, we first project them onto a reduced latent space via the symplectic encoder $\fenc:\mathbb{R}^{2n}\rightarrow\mathbb{R}^{2k}$:
\begin{equation}
\mathbf{y}_i = \fenc(\mathbf{x}_i)\in \mathbb{R}^{2k}, \quad i=0,\dots,T.
\end{equation}

Subsequently, the temporal evolution within this latent space is approximated by an autonomous, discrete-time symplectic flow map $\fflow:\mathbb{R}^{2k}\to\mathbb{R}^{2k}$. We primarily construct $\fflow$ using a latent-space H\'enonNet, denoted as $\mathcal{H}_{\mathrm{flow}}$:
\begin{equation}
\fflow = \mathcal{H}_{\mathrm{flow}}.
\end{equation}
The discrete evolution in latent coordinates is then represented as:
\begin{equation}
    \mathbf{y}_{i+1} = \fflow(\mathbf{y}_i).
\end{equation}


We train this latent dynamics module by minimizing the flow prediction error:
\begin{equation}
    \mathcal{L}_f = \frac{1}{T}\sum_{i=1}^{T}\|\mathbf{y}_i - \fflow(\mathbf{y}_{i-1})\|_2^2.
\end{equation}
Due to the intrinsic symplectic structure of the latent-space H\'enonNet ($\mathcal{H}_{\mathrm{flow}}$), the trained flow map $\fflow$ inherently preserves the symplectic structure of the reduced-order latent system. 


\subsection{Unified training process}
Although the autoencoder and flow map can be trained separately \cite{simpson2021machine,reddy2019reduced}, we adopt an end-to-end training strategy that jointly optimizes dimensionality reduction and latent dynamics modeling. This unified framework ensures dynamic consistency and physical interpretability within the learned latent representations. Specifically, we introduce the following combined loss function to concurrently promote reconstruction fidelity and predictive accuracy:



\begin{equation}\label{eq:loss-rom}
\mathcal{L}_{\mathrm{rom}}
=
\sum_{i=0}^{T}
\left\|\mathbf{x}_i-\fdec(\fenc(\mathbf{x}_i))\right\|_2
+
\sum_{i=0}^{T-1}
\left\|\fdec\bigl(\fflow(\fenc(\mathbf{x}_i))\bigr)-\mathbf{x}_{i+1}\right\|_2 .
\end{equation}
While our architecture inherently ensures symplecticity, exact Hamiltonian conservation is not enforced, as the proposed architecture only guarantees an approximate Hamiltonian is enforced. To explicitly maintain Hamiltonian invariance, particularly crucial for complex nonlinear systems, we introduce an additional Hamiltonian conservation constraint:

\begin{equation}\label{eq:loss-ham}
\begin{split}
     \mathcal{L}_{\rm ham} = & \sum_{i=0}^{T} \left\|\tilde{H}(\fdec( \fenc(\mathbf{x}_i)))-\mathit{H}(\mathbf{x}_0)\right\|_2 \\
    &  \quad + \sum_{i=0}^{T-1}\left\|\mathit{H}(\fdec( \fflow( \fenc(\mathbf{x}_i))))-\tilde{H}(\mathbf{x_0})\right\|_2
\end{split}
\end{equation}
where $\tilde{H}$ denotes the original Hamiltonian of the system state $\mathbf{x}$.
This constraint ensures Hamiltonian conservation during both reconstruction and latent space dynamics evolution.
The complete loss function is therefore expressed as:
\begin{equation}\label{eq:loss_total}
     \mathcal{L}_{\rm total} = \lambda_1\mathcal{L}_{\rm rom} + \lambda_2 \mathcal{L}_{\rm ham}
\end{equation}
where $\lambda_1$ and $\lambda_2$ are hyperparameters that balance the ROM loss $\mathcal{L}_{\rm rom}$ and the Hamiltonian conservation loss $\mathcal{L}_{\rm ham}$.  To further enhance stability and generalization, we incorporate two training strategies:

\paragraph{Multi-step training strategy} The basic one-step training approach predicts the state at time step $i+1$ using the true value at time step $i$:
\begin{equation}
    \hat{\mathbf{x}}_{i+1} = \fdec(\fflow(\fenc({\mathbf{x}}_{i})))
\end{equation}

Although computationally efficient, this approach may exhibit limited generalization capability. To enhance prediction accuracy and stability, we implement an auto-regressive training scheme that unfolds over multiple time steps. This approach uses the model's prediction $\hat{\mathbf{x}}_{i}$ at time step $i$ as input for the subsequent time step prediction:
\begin{equation}
    \hat{\mathbf{x}}_{i+1} = \fdec(\fflow(\fenc(\hat{\mathbf{x}}_{i})))
\end{equation}

For an $M$-step unrolling during training, the ROM loss \eqref{eq:loss-rom} is reformulated as:

\begin{equation}\label{eq:loss-rom-new}
\begin{split}
\mathcal{L}_{\mathrm{rom}}
={}&
\sum_{i=0}^{T}\left\|\mathbf{x}_i-\fdec\circ \fenc(\mathbf{x}_i)\right\|_2 \\
&\quad
+\sum_{i=0}^{T-M-1}\sum_{j=0}^{M-1}
\left\|\fdec\circ \fflow\circ \fenc(\hat{\mathbf{x}}_{i+j})-\mathbf{x}_{i+j+1}\right\|_2 .
\end{split}
\end{equation}
while the Hamiltonian conservation loss takes the form:

\begin{equation}\label{eq:loss-ham-new}
\begin{split}
\mathcal{L}_{\mathrm{ham}}
={}&
\sum_{i=0}^{T}\left\|\mathit{H}\left(\fdec\circ \fenc(\mathbf{x}_i)\right)-\mathit{H}(\mathbf{x}_0)\right\|_2 \\
&\quad
+\sum_{i=0}^{T-M-1}\sum_{j=0}^{M-1}
\left\|\mathit{H}\left(\fdec\circ \fflow\circ \fenc(\hat{\mathbf{x}}_{i+j})\right)-\mathit{H}(\mathbf{x}_0)\right\|_2 .
\end{split}
\end{equation}
Here, the state evolution $\hat{\mathbf{x}}_{i+j}$ follows a recursive relationship defined as:
\begin{equation}
    \hat{\mathbf{x}}_{i+j} = 
    \begin{cases}
        \mathbf{x}_{i+j}, & j=0, \\
        \fdec(\fflow(\fenc(\hat{\mathbf{x}}_{i+j-1}))), & j>0.
    \end{cases}
\end{equation}

Although this strategy increases training complexity, it significantly enhances prediction stability, as demonstrated in \cite{brandstetter2021message}.


\paragraph{Noise injection strategy} Small-amplitude Gaussian perturbations $\epsilon\sim\mathcal{N}(0,\sigma^2)$ are added to the inputs during training to improve robustness and balance optimization between competing objectives \cite{citation-key,pfafflearning,sanchez2020learning} as follows,
\begin{equation}
    \tilde{\mathbf{x}}_{i} = \mathbf{x}_{i} + \epsilon,\quad i=0,\dots,T-M-1.
\end{equation}



\section{Numerical results} \label{sec:results}

In this section, we demonstrate the effectiveness of our symplectic ROM framework through three canonical Hamiltonian problems: the linear wave equation, parametric linear wave equation, and nonlinear Schr{\"o}dinger equation. These examples are chosen to validate our approach across different levels of complexity, from linear dynamics to parametric dependence and nonlinear behavior.

Our symplectic model comprises two main components: a symplectic autoencoder for dimensional reduction and a latent flow mapping module for time evolution. The autoencoder consists of H\'enonNet and optional G-reflectors, while the latent flow mapping employs a similar architecture but with reduced dimensionality. For each test problem, we provide details of the architecture in the corresponding network structure tables.

In all numerical examples, an important metric we are interested in investigating is the approximation power of the proposed symplectic autoencoders. 
We focus on four types of parametrization: the cotangent lift (embedding) proposed in \cite{peng2016}, the G-reflector-enhanced linear embedding proposed in Section~\ref{sec:linear},
the H\'enonNet-based embedding proposed in Section~\ref{sec:nonlinear},
and the symplectic embedding based on the composition of G-reflectors and H\'enonNets discussed in Section~\ref{sec:composition}.
It is worth emphasizing again that the cotangent lift is an incomplete parameterization of $\mathrm{SpSt}(2k,2n)$, whereas the others offer complete parameterizations of the corresponding linear or non-linear symplectic embedding manifolds.
To define a unified metric, given a snapshot sequence $\{\mathbf{x}_i\}_{i=0}^N$ and its corresponding reconstructions $\{\hat{\mathbf{x}}_i\}_{i=0}^N$,
we measure the reconstruction loss by the mean squared error (MSE):
\begin{equation}
  \mathrm{MSE}
  = \frac{1}{N+1} \sum_{i=0}^N \bigl\|\mathbf{x}_i - \hat{\mathbf{x}}_i\bigr\|_2^2.
\end{equation}

\paragraph{Implementation details}
All base neural networks use a multilayer perceptron (MLP) and a exponential linear unit (ELU) \cite{clevert2016fast} activation function. The H\'enonNet utilizes fully connected layers with the specified width and depth, while G-reflectors only needs to specify the layers. Training is performed using the Adam optimizer \cite{Ilya_fix_2017} with an initial learning rate of $10^{-3}$ and exponential decay throughout training. For reference data generation, we employ the symplectic  St{\"o}rmer-Verlet integration scheme. To assess method robustness under realistic conditions, we inject additive Gaussian noise into training data with magnitude $\epsilon \sim 10^{-3}\times \mathcal{N}(0,1)$, representing typical measurement uncertainties while preserving solution fidelity. 
Throughout this study, we present comprehensive results using the symplectic neural network approach (denoted as ``Sym''), while providing architectural comparisons specifically for reconstruction accuracy assessment. All experimental results compare our symplectic method against high-fidelity reference solutions computed using the St{\"o}rmer-Verlet scheme (denoted as ``Exact''), ensuring rigorous validation against established symplectic integration standards. For all trainable neural architectures, evaluation is performed on independently generated test samples, and the reported test reconstruction errors are averaged over three runs with different random seeds.

\begin{example}\label{ex:linear-wave} We first consider the following  1D linear wave equation
\begin{equation}
   \left\{\begin{array}{ll}
u_{t t}(x, t)=\omega^2 u_{x x}(x, t), \quad x\in\Omega=[0,1],\\
u(x, 0)=u^{0}(x),\\
u_t(x,0) = -\omega u^{0}_x(x),\\
u(0, t)=u(1, t)=0.
\end{array}\right. 
\end{equation}
The initial condition is constructed using the cubic spline function $h(s)$ defined as
\begin{equation}
    h(s)=\left\{\begin{array}{lr}
1-\frac{3}{2} s^{2}+\frac{3}{4} s^{3}, & 0 \leq s \leq 1 \\
\frac{1}{4}(2-s)^{3}, & 1<s \leq 2 \\
0, & s>2
\end{array}\right.\label{eqn:linear1}
\end{equation}
with the initial condition given by
\begin{equation}
    u^{0}(x) = h\left( a_0 \left|x - \frac12-x_0\right| \right).
\end{equation}
where $a_0\in [7,9]$, $x_0\in[-0.2,0.2]$, and $\omega$ represents the wave speed. 
The exact solution to \eqref{eqn:linear1} is 
\begin{equation}
    u(x,t
) =u^{0}(x-\omega t).
\end{equation}

To facilitate a symplectic formulation, we introduce new variables $q = u$ and $p = q_t$, leading to the first-order system
\begin{equation}
    \left\{\begin{array}{l}
q_{t}(x, t, \omega)=p(x, t, \omega) \\
p_{t}(x, t, \omega)=\omega^2 q_{x x}(x, t, \omega)
\end{array}\right.
\end{equation}
with the associated continuous Hamiltonian given by
\begin{equation}
    H_{\rm cont}(q, p)=\frac{1}{2} \int_{0}^{1} p^{2}+\omega^2 q_{x}^{2} d x.
\end{equation}

For numerical discretization, the spatial domain is divided uniformly into $N$ grid points $x_i=i\Delta x$, with mesh size $\Delta x=1/N$. Denoting $q_i = q(t, x_i, \omega)$, $p_i = p(t, x_i, \omega)$, and introducing the state vector $\mathbf{z} = (q_1,\ldots,q_N,p_1,\ldots,p_N)$, we obtain the discrete system through second-order central finite difference,
\begin{equation}
    \frac{d}{d t} \mathbf{z}=\bbJ_{2 N} \nabla_{\bf z} H({\bf z})=\bbJ_{2 N} L \, \mathbf{z}
\end{equation}
where the discrete Hamiltonian is given by
\begin{equation}
    H(\mathbf{z}) = \frac{\Delta x}{2} \sum_{i=1}^{N}\left(p_{i}^{2}+ \omega^2 \frac{\left(q_{i+1}-q_{i}\right)^{2}}{2 \Delta x^{2}}+ \omega^2 \frac{\left(q_{i}-q_{i-1}\right)^{2}}{2 \Delta x^{2}}\right),
\end{equation}
and
 \begin{equation}
     L=\left(\begin{array}{cc}
I_{N} & 0_{N} \\
0_{N} & \omega^2 D_{x x}
\end{array}\right).
 \end{equation}
 Here the difference operator $D_{xx}$ is defined as: 
 \begin{equation}\label{eq:D_xx}
D_{xx}
=\frac{1}{\Delta x^{2}}
\begin{pmatrix}
-2 & 1  &        &        &        \\
 1 & -2 & 1      &        &        \\
   & \ddots & \ddots & \ddots &   \\
   &        & 1      & -2     & 1  \\
   &        &        & 1      & -2
\end{pmatrix}_{\!N\times N}.
\end{equation}

\begin{table}[H]
\centering
\caption{Numerical discretization parameters used in Example~\ref{ex:linear-wave}.}
\label{table:ex5-1-config}
\small
\begin{tabular}{|l|c|}
\hline
Space discretization size & $\Delta x = 0.005$ \\ \hline
Time discretization size  & $\Delta t = 0.24$ \\ \hline
Wave speed                & $\omega^2 = 0.01$ \\ \hline
\end{tabular}
\end{table}

\begin{table}[H]
\centering
\caption{Neural network architectures used in Example~\ref{ex:linear-wave}.}
\label{table:ex5-1-network}
\small
\begin{tabular}{|l|c|}
\hline
Autoencoder H\'enonNet      & $[300]\times 2$ \\ \hline
Autoencoder G-reflector     & $10$ \\ \hline
Flow mapping H\'enonNet     & $[30]\times 2$ \\ \hline
\end{tabular}
\end{table}

Numerical experiments utilize the parameters provided in Table \ref{table:ex5-1-config}. We generate snapshot data by uniformly sampling 3000 parameter pairs $(a_0,x_0)$ from $[7,9]\times[-0.2,0.2]$ and integrating numerically over $[0,12]$ using the St{\"o}rmer-Verlet scheme. In addition, we sample 300 parameter pairs from the same parameter range as an independent test set for reconstruction evaluation. The latent dimension is set to 10, and neural network architectures are summarized in Table \ref{table:ex5-1-network}. Hyperparameters $(\lambda_1,\lambda_2)$ in the loss function \eqref{eq:loss_total} are $(1,0.01)$. 

\begin{figure}[!htbp]
 \centering
     \subfigure[Latent space dynamics]{\includegraphics[width=0.45\textwidth]{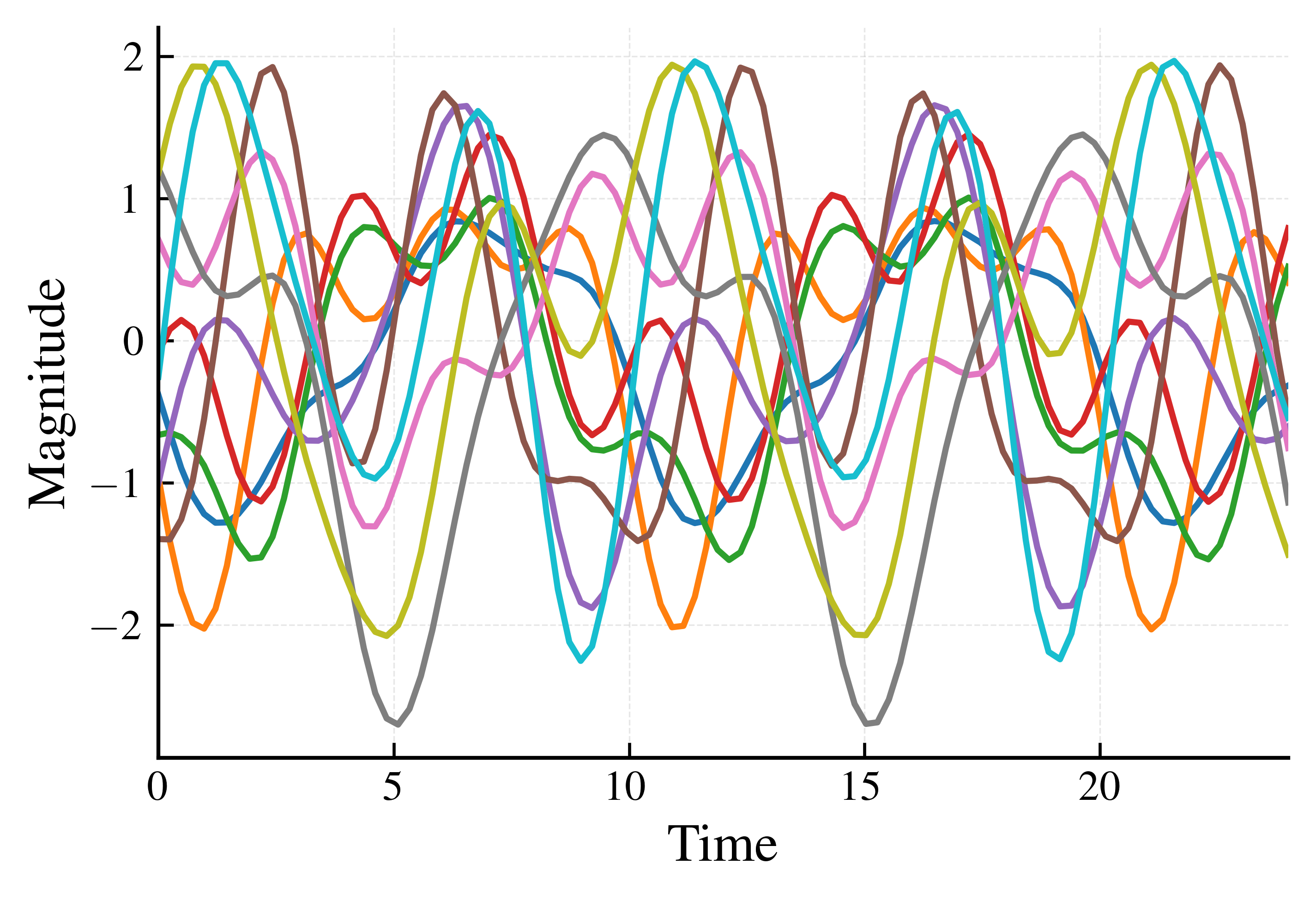}}\quad
     \subfigure[$t=0$]{\includegraphics[width=0.45\textwidth,]{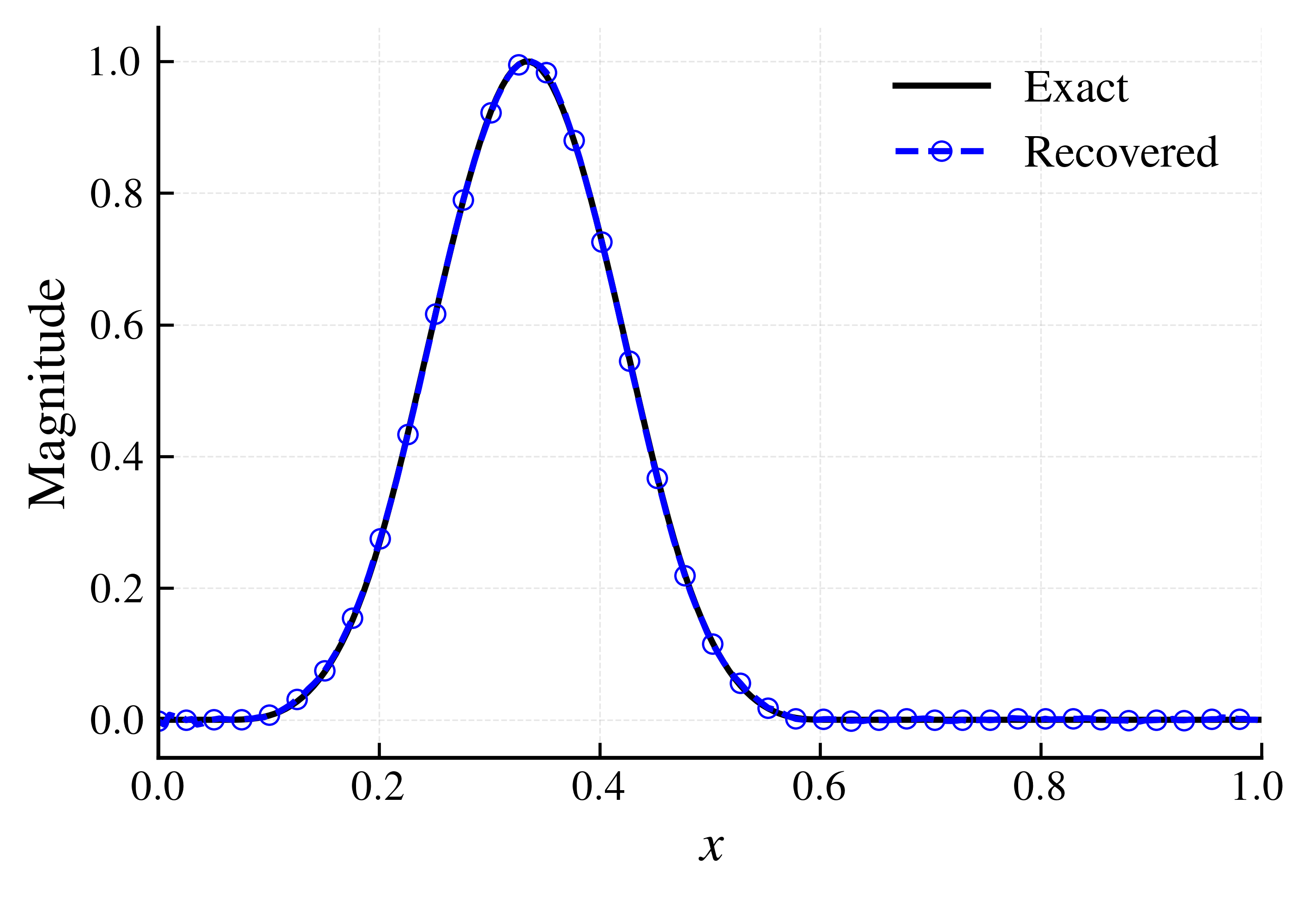}}\\
     \subfigure[$t=12$]{\includegraphics[width=0.45\textwidth]{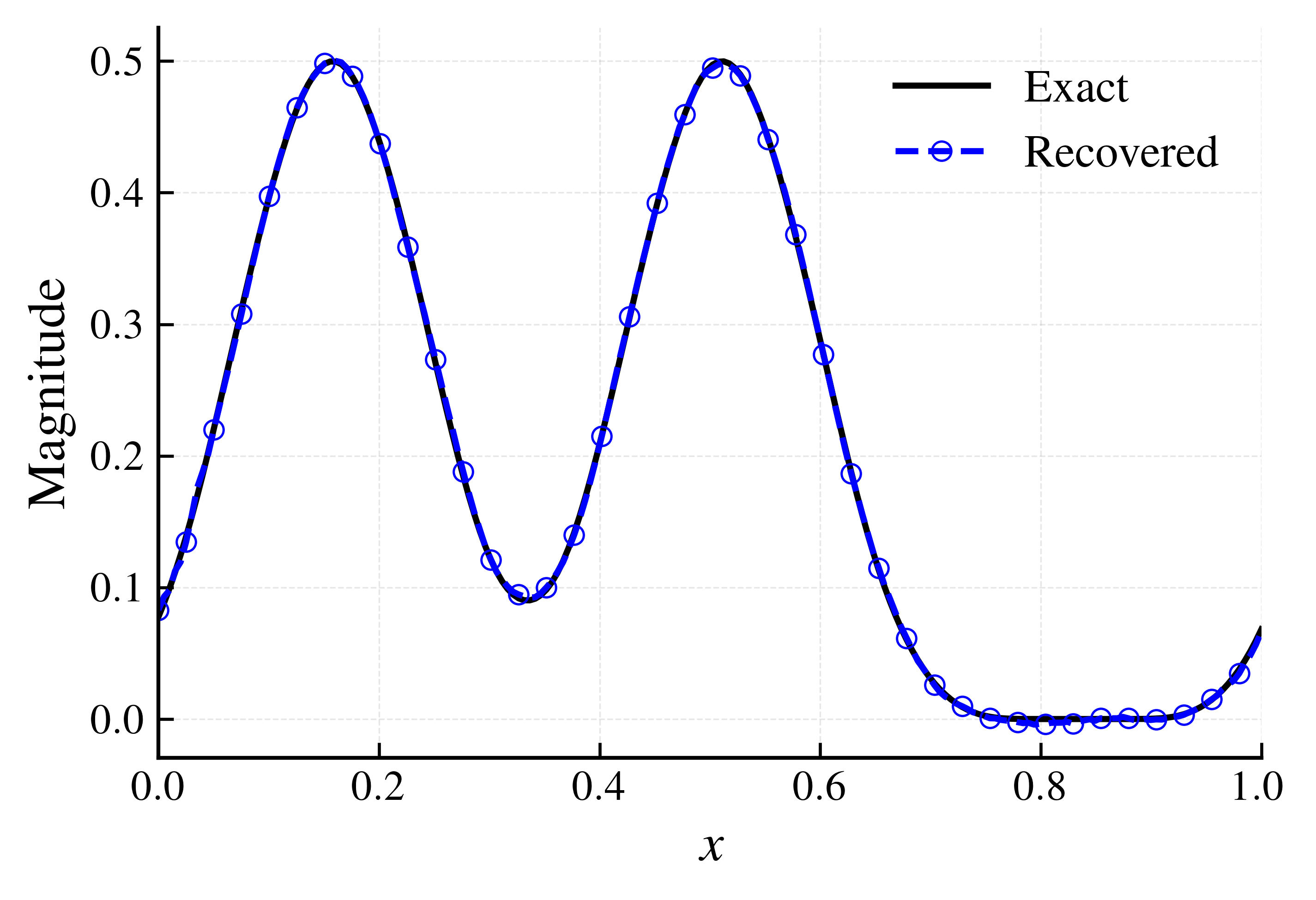}}\quad
     \subfigure[$t=24$]{\includegraphics[width=0.45\textwidth]{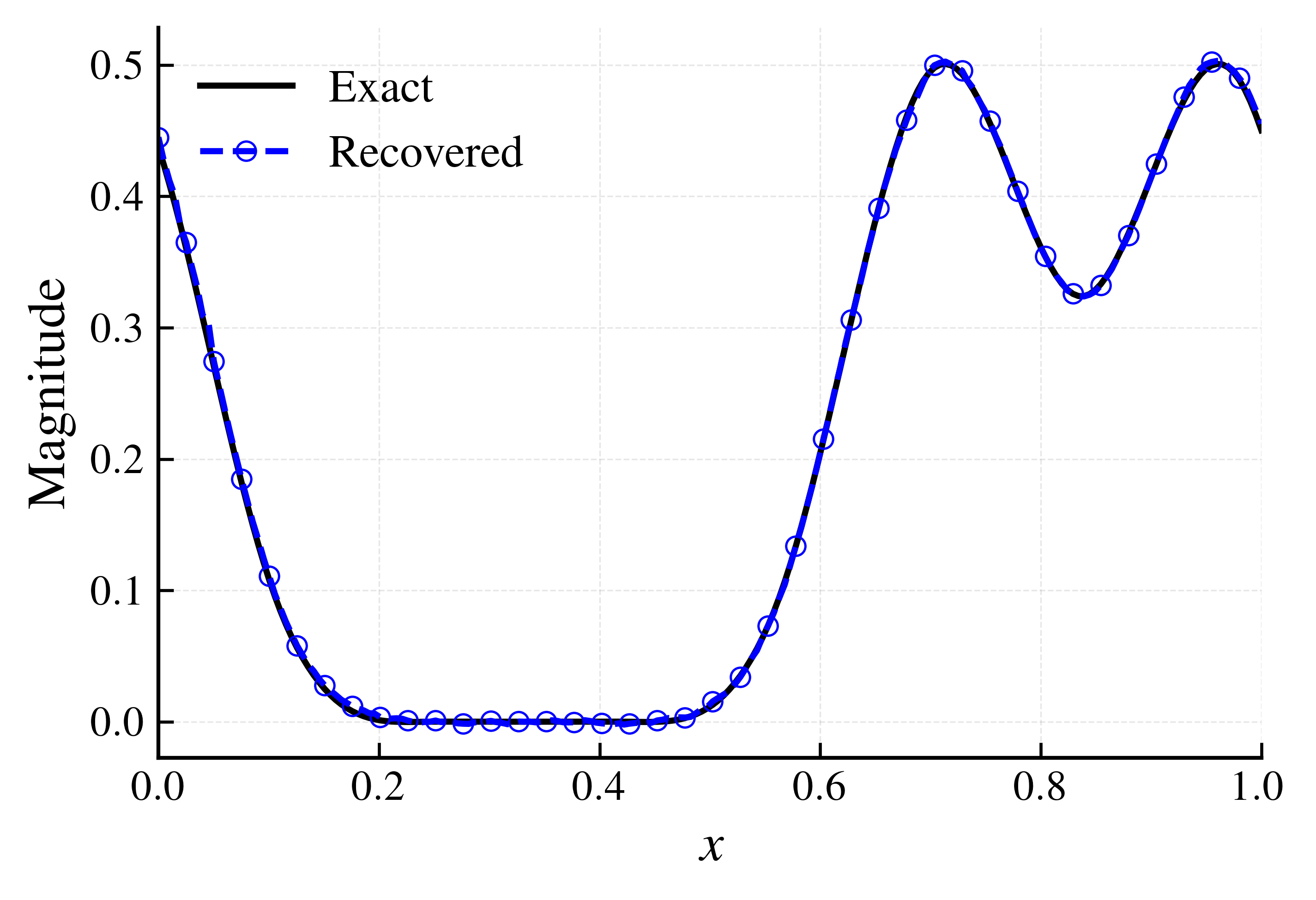}}
	\caption{ROM prediction in Example~\ref{ex:linear-wave} for $(a_0,x_0)=(8.939,0.079)$: (a) latent space dynamics; (b)--(d) reconstructed wave solutions at $t=0$, $t=12$, and $t=24$.}\label{fig:linear-1}
\end{figure}

Figure~\ref{fig:linear-1} demonstrates the performance of our method on a representative case with $a_0 = 8.939$ and $x_0 = 0.079$. 
Figure \ref{fig:linear-1}(a) shows that the latent space trajectories exhibit smooth, stable evolution, while Figures \ref{fig:linear-1}(b-d) show that the reconstructed solutions are stable and accurate across multiple time instances. Notably, despite training only on data up to $t = 12$, the model successfully extrapolates to $t = 24$, demonstrating robust long-term prediction capability.
The conservation properties of our symplectic method are observed in Figure \ref{fig:linear-2}, which displays the time evolution of the discrete Hamiltonian along a solution trajectory and the instantaneous Hamiltonian deviation from its initial value. The result confirms that our method preserves the Hamiltonian structure of the system. 

\begin{table}[H]
\centering
\caption{Test reconstruction MSE in Example~\ref{ex:linear-wave}, including Hamiltonian Operator Inference (H-OpInf) baseline. For trainable neural architectures, the reported values are mean $\pm$ standard deviation over three independent runs; the cotangent-lift and H-OpInf results are deterministic.}
\label{tab:linear-rec}
\small
\setlength{\tabcolsep}{3pt}  

\begin{tabular}{|l|c|c|c|c|c|}
\hline
Method & Cotangent-lift & G-reflector & H\'enonNet & H\'enonNet + G-reflector & H-OpInf \\ \hline
MSE & 3.260E-4 & 
\makecell{2.759E-4 \\ $\pm$ \scriptsize{1.248E-4}} & 
\makecell{7.292E-6 \\ $\pm$ \scriptsize{7.644E-7}} & 
\makecell{3.470E-6 \\ $\pm$ \scriptsize{6.141E-7}} & 
8.368E-5 \\ \hline
\end{tabular}
\end{table}

  Furthermore, to provide a comprehensive assessment of reconstruction accuracy, we compare five distinct approaches: the cotangent-lift method, a 100-layer G-reflector implementation, a standalone H\'enonNet architecture, a composite H\'enonNet + G-reflector architecture, and the Hamiltonian Operator Inference (H-OpInf) baseline \cite{gruber2025variationally}. For the H-OpInf baseline, we examined several H-OpInf-style basis families and report the best reconstruction performance achieved using a centered POD basis. Table~\ref{tab:linear-rec} summarizes the resulting test reconstruction MSEs. The results exhibit a clear performance hierarchy: although the standalone H\'enonNet already yields a substantial improvement over the linear baselines, the combined H\'enonNet + G-reflector architecture attains the lowest reconstruction error. By contrast, the H-OpInf baseline remains less accurate than all symplectic ROM variants considered here. For this linear wave problem, where strong linear structure is present, the additional linear symplectic corrections introduced by the G-reflector further enhance the effectiveness of the H\'enonNet-based nonlinear representation.

\begin{figure}[!htbp]
 \centering
     \subfigure[Hamiltonian]{\includegraphics[width=0.45\textwidth]{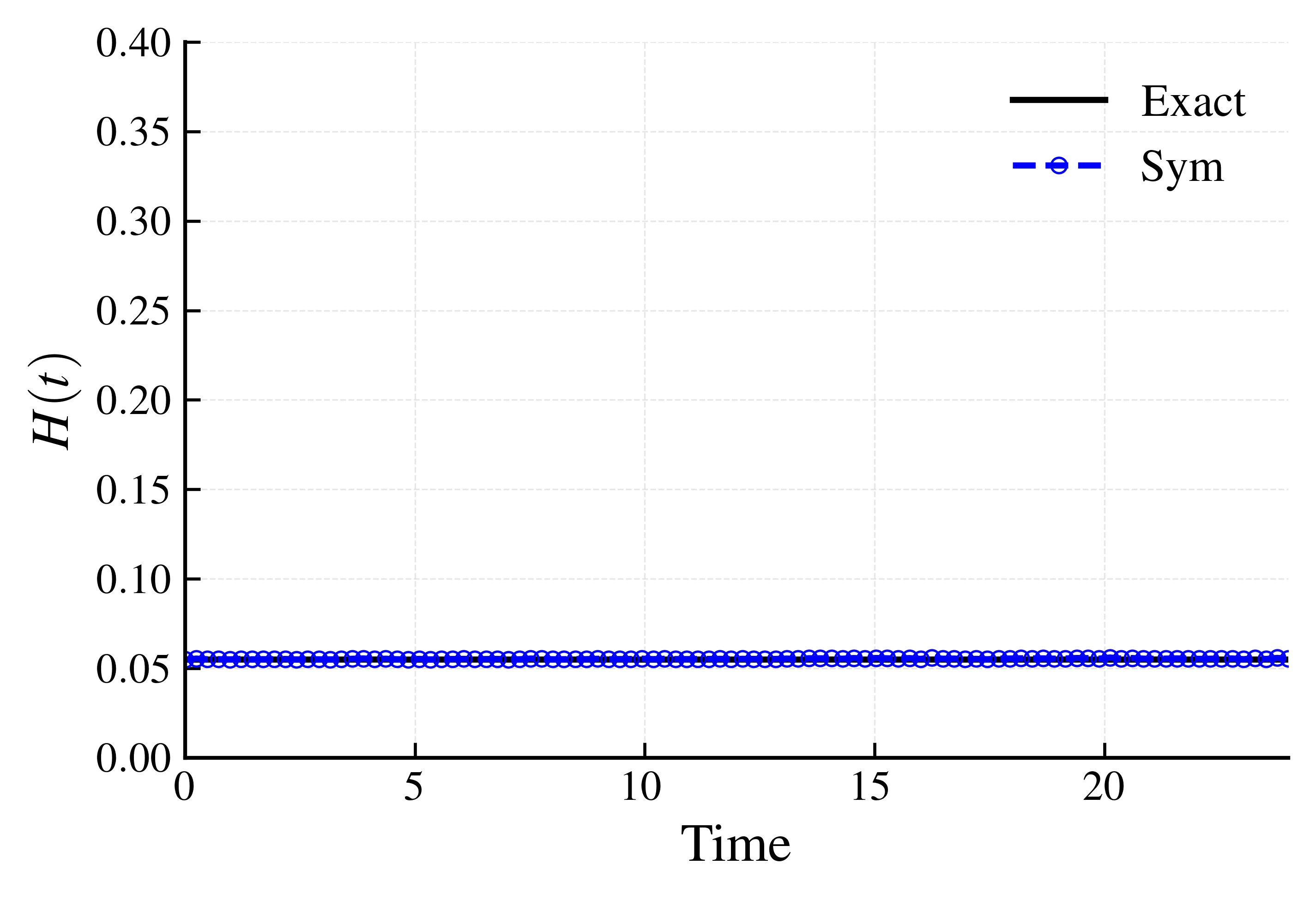}}\quad\quad
     \subfigure[Hamiltonian deviation]{\includegraphics[width=0.45\textwidth]{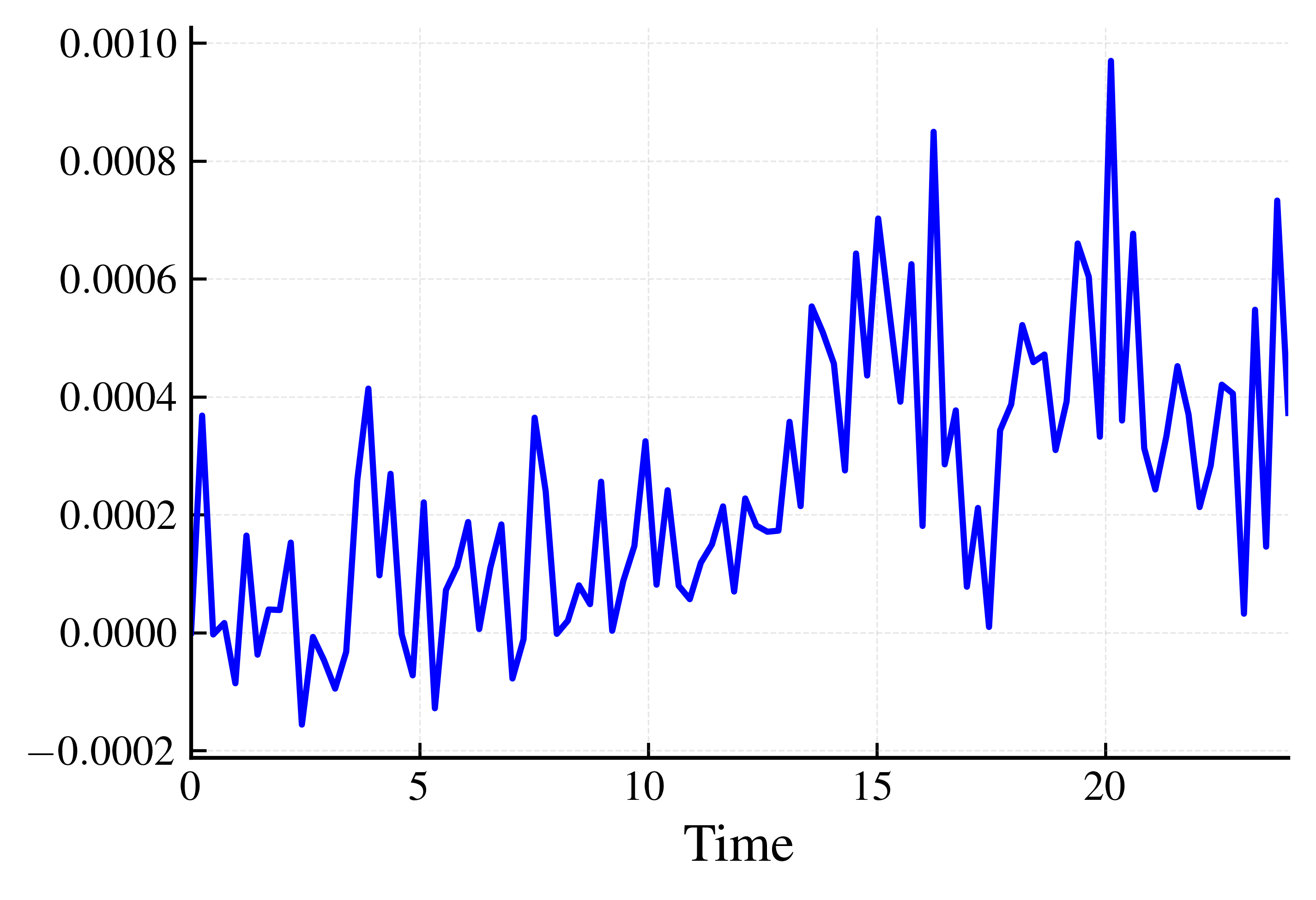}}
	\caption{Hamiltonian conservation results in Example~\ref{ex:linear-wave}: (a) discrete Hamiltonian in time; (b) instantaneous Hamiltonian deviation from the initial value.}\label{fig:linear-2}
\end{figure} 

\end{example}

\begin{example}\label{ex:param-linear}
We next consider the parametric linear wave equation:
\begin{equation}
    \left\{\begin{array}{l}
u_{t t}(x, t, \omega)=\kappa(\omega) u_{x x}(x, t, \omega), \quad x\in[0,1],\\
u(x,0)= u^0(x),\\
\end{array}\right.
\end{equation}
where the wave speed parameter $\kappa(\omega)$ depends on a four-dimensional parameter vector $\omega=(\omega_1,\omega_2,\omega_3,\omega_4)$ via:
\[\kappa(\omega) = c^2\left(\sum_{l=1}^{4}\frac{1}{l^2}\omega_l\right),\]
and $c\in\mathbb{R}$ is a fixed constant.

The initial condition is defined through the cubic spline function $h(s)$ with
\begin{equation}
    h(s)=\left\{\begin{array}{lr}
1-\frac{3}{2} s^{2}+\frac{3}{4} s^{3}, & 0 \leq s \leq 1 \\
\frac{1}{4}(2-s)^{3}, & 1<s \leq 2 \\
0, & s>2
\end{array}\right.
\end{equation}
and given by
\begin{equation}
    u^{0}(x) = h\left( 10 \left|x -\frac12\right| \right).
\end{equation}

Introducing canonical variables $q=u$ and $p=q_t$, the system can be rewritten in the canonical form
\begin{equation}
    \left\{\begin{array}{l}
q_{t}(x, t, \omega)=p(x, t, \omega) \\
p_{t}(x, t, \omega)=\kappa(\omega) q_{x x}(x, t, \omega)
\end{array}\right.
\end{equation}
with the continuous Hamiltonian given by
\begin{equation}
    H_{\rm cont}(q, p, \omega)=\frac{1}{2} \int_{0}^{1} p^{2}+\kappa(\omega) q_{x}^{2} d x.
\end{equation}

The numerical discretization employs a uniform spatial grid with $N$ points, $x_i=i\Delta x$, $\Delta x=1/N$. Defining discrete variables $q_i=q(t,x_i,\omega)$, $p_i=p(t,x_i,\omega)$, and the state vector $\mathbf{z}=(q_1,\dots,q_N,p_1,\dots,p_N)$, the discretized system obtained via second-order central finite differences is
\begin{equation}
    \frac{d}{d t} \mathbf{z}=\bbJ_{2 N} \nabla_{\bf z} H({\bf z})=\bbJ_{2 N} L\mathbf{z},
\end{equation}
where the discrete Hamiltonian is given by
\begin{equation}
    H(\mathbf{z}) = \frac{\Delta x}{2} \sum_{i=1}^{N}\left(p_{i}^{2}+\kappa(\omega) \frac{\left(q_{i+1}-q_{i}\right)^{2}}{2 \Delta x^{2}}+\kappa(\omega) \frac{\left(q_{i}-q_{i-1}\right)^{2}}{2 \Delta x^{2}}\right)
\end{equation}
and matrix $L$ is
 \begin{equation}
     L=\left(\begin{array}{cc}
I_{N} & 0_{N} \\
0_{N} & \kappa(\omega) D_{x x}
\end{array}\right).
 \end{equation}

\begin{table}[H]
\centering
\caption{Numerical discretization parameters used in Example~\ref{ex:param-linear}.}
\label{table:ex5-2-config}
\small
\begin{tabular}{|l|c|}
\hline
Space discretization size & $\Delta x = 0.005$ \\ \hline
Time discretization size  & $\Delta t = 0.24$ \\ \hline
Wave speed                & $c^2 = 0.1$ \\ \hline
\end{tabular}
\end{table}

\begin{table}[H]
\centering
\caption{Neural network architectures used in Example~\ref{ex:param-linear}.}
\label{table:ex5-2-network}
\small
\begin{tabular}{|l|c|}
\hline
Autoencoder H\'enonNet      & $[512]\times 2$ \\ \hline
Autoencoder G-reflector     & $10$ \\ \hline
Flow mapping H\'enonNet     & $[32]\times 2$ \\ \hline
\end{tabular}
\end{table}

Numerical experiments use the discretization parameters given in Table~\ref{table:ex5-2-config}. Snapshot data are generated by sampling the four-dimensional parameter space $[0,1]^4$ uniformly using a grid of $5^4$ points. Numerical trajectories are computed over the interval $[0,3]$ using the St{\"o}rmer--Verlet scheme. In addition, we randomly sample 100 parameter instances from the same parameter space as an independent test set for reconstruction evaluation. The latent dimension is 10, with the neural network architecture detailed in Table~\ref{table:ex5-2-network}. The hyperparameters $\lambda_1$ and $\lambda_2$ in the total loss function \eqref{eq:loss_total} are assigned values of $1$ and $0.01$, respectively. 


\begin{figure}[!htbp]
 \centering
     \subfigure[Latent space dynamics]{\includegraphics[width=0.45\textwidth]{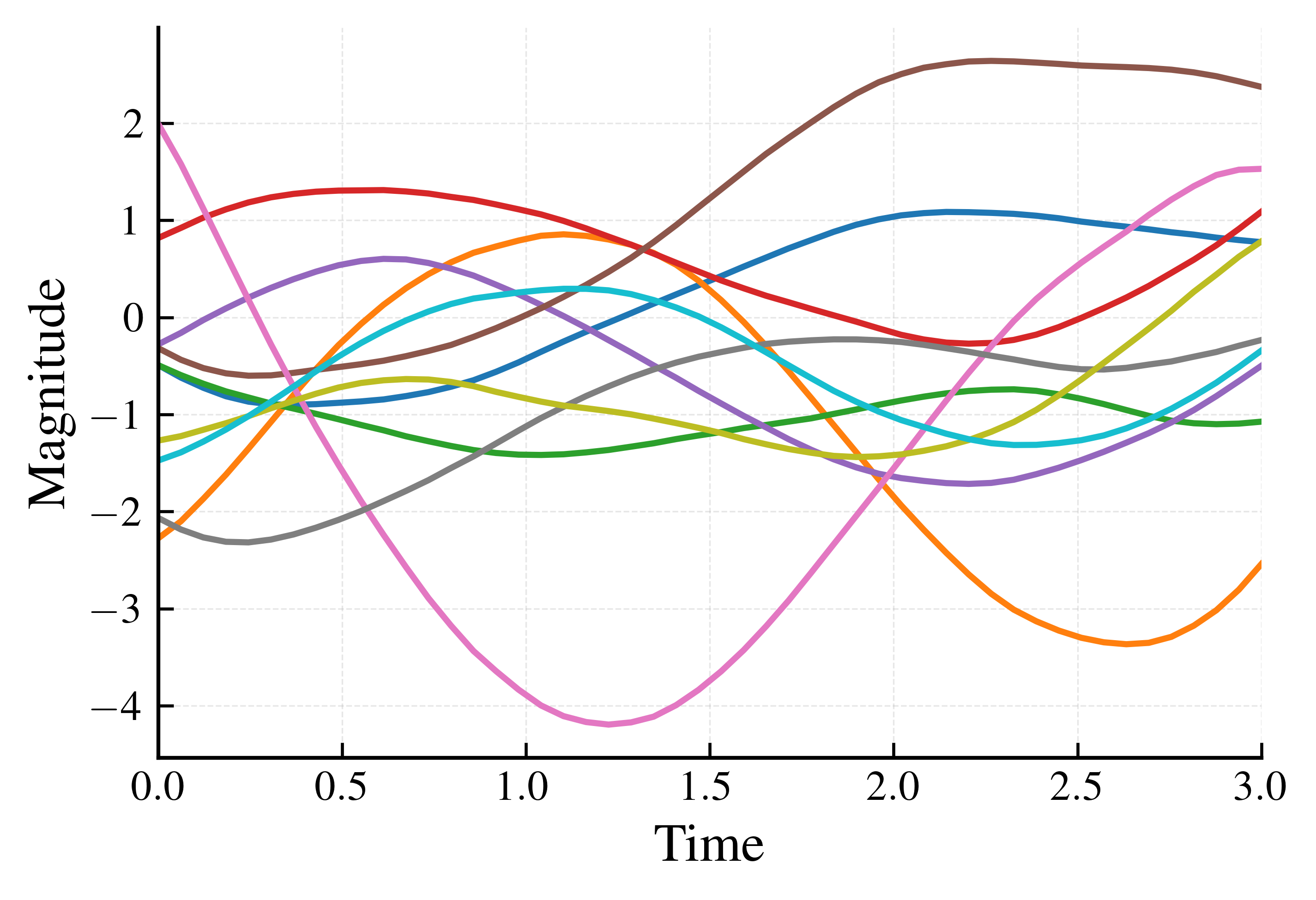}}\quad
     \subfigure[$t=0$]{\includegraphics[width=0.45\textwidth]{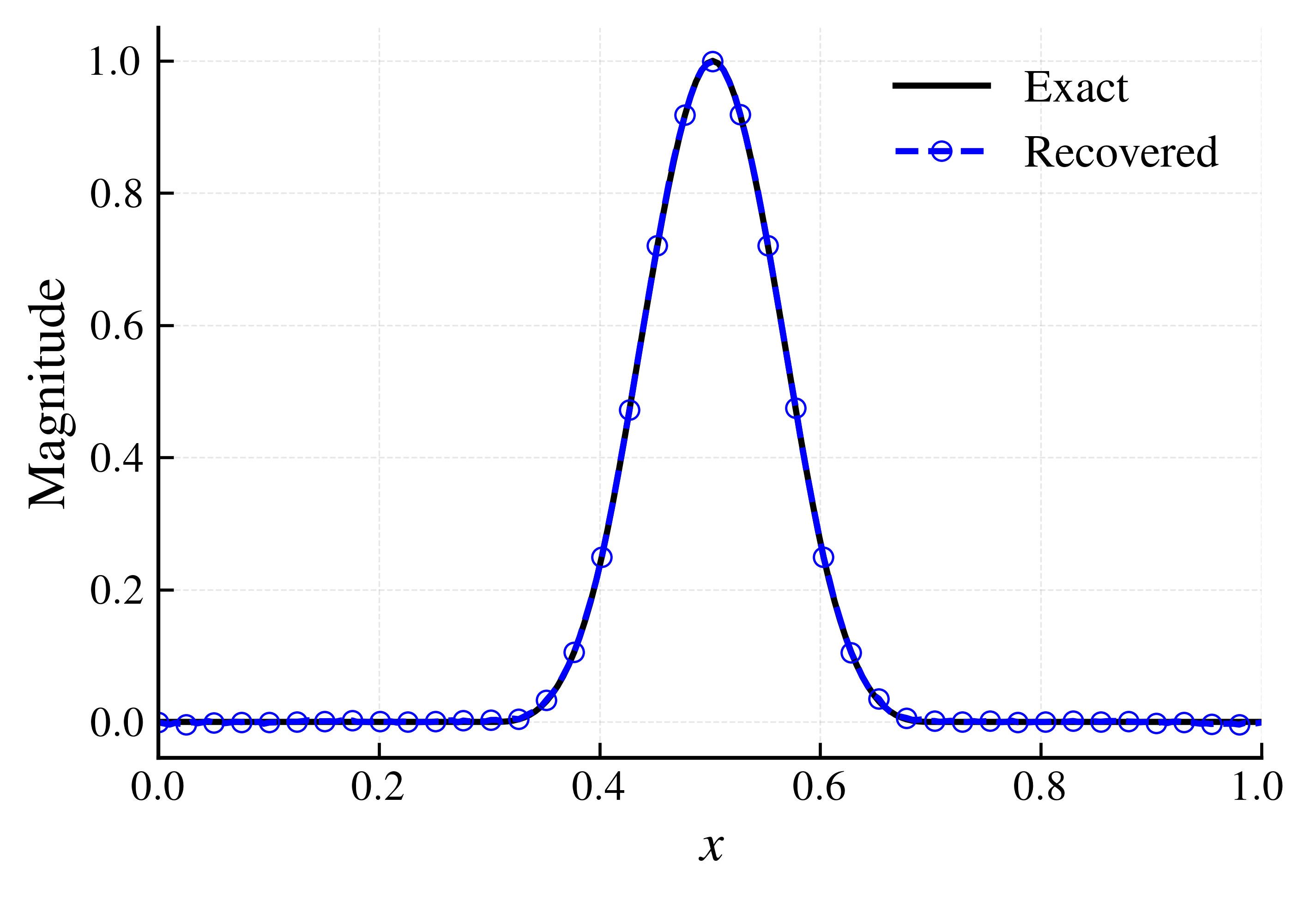}}\\
     \subfigure[$t=1.5$]{\includegraphics[width=0.45\textwidth]{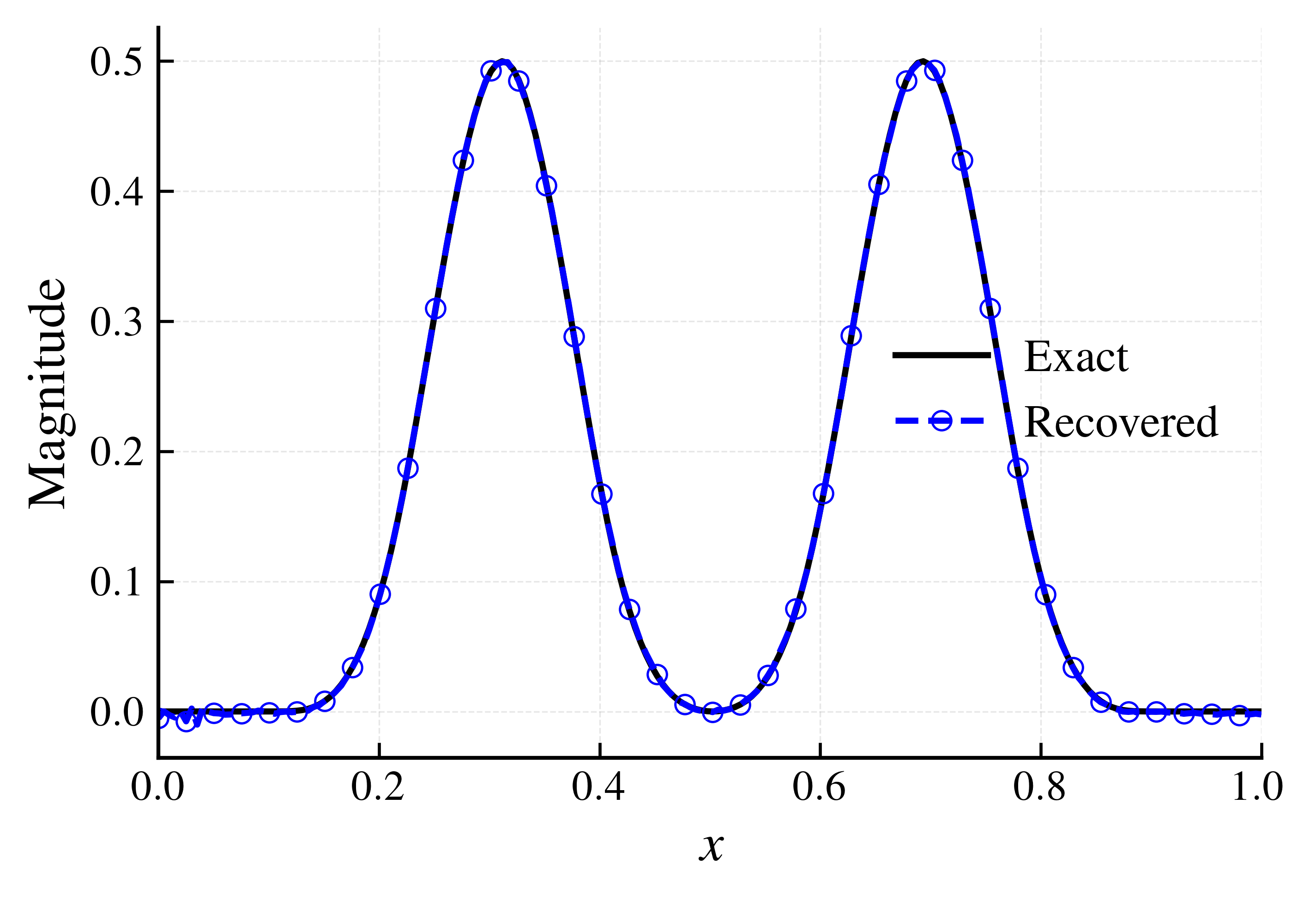}}\quad
     \subfigure[$t=3$]{\includegraphics[width=0.45\textwidth]{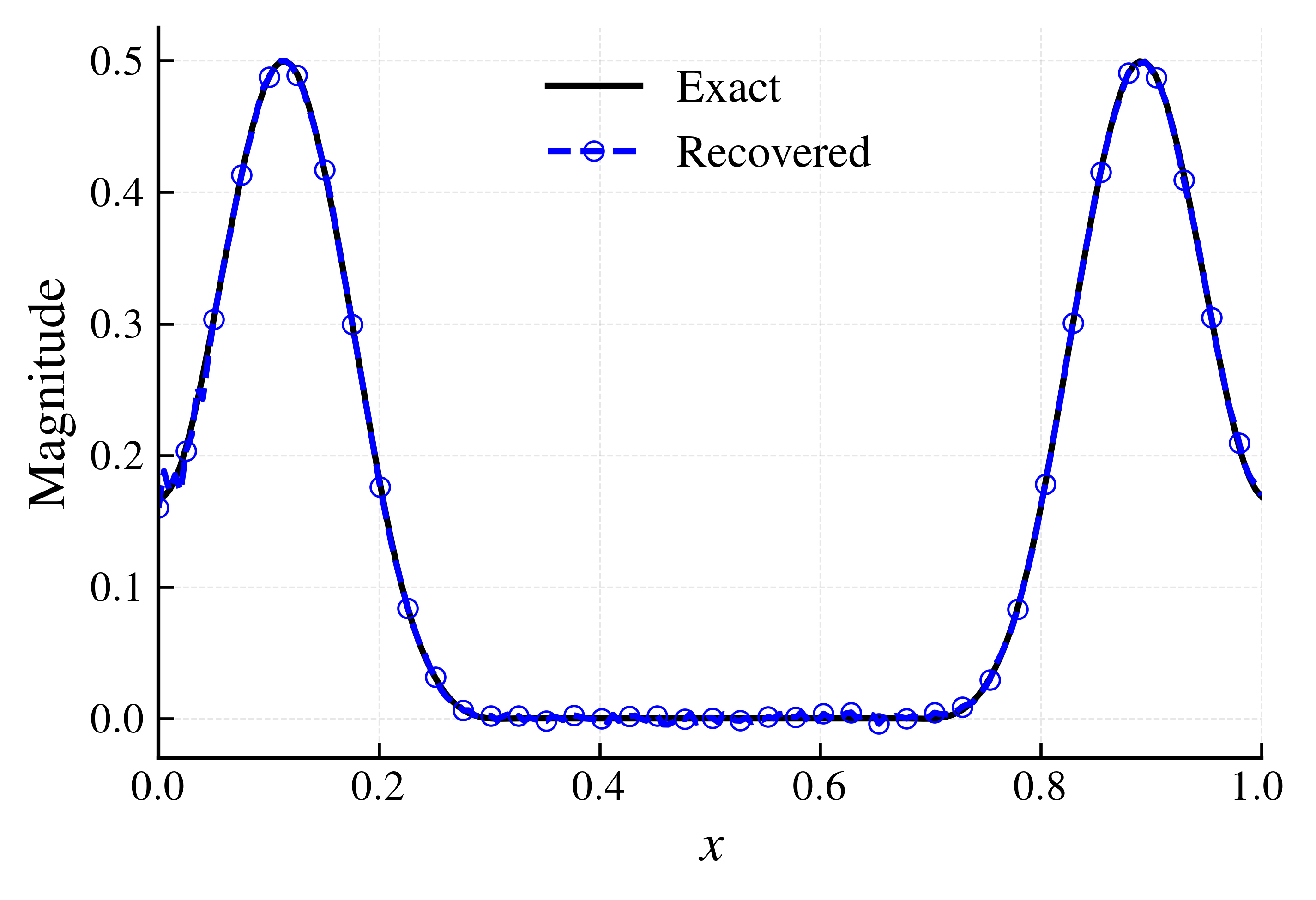}}
	\caption{ROM prediction in  Example~\ref{ex:param-linear} with $\kappa(\omega)=0.01736$: (a) latent space dynamics; reconstructed solutions at (b) $t=0$, (c) $t=1.5$, and (d) $t=3$.}\label{fig:param_linear}
\end{figure}

\begin{figure}[!htbp]
 \centering
     \subfigure[Hamiltonian]{\includegraphics[width=0.45\textwidth]{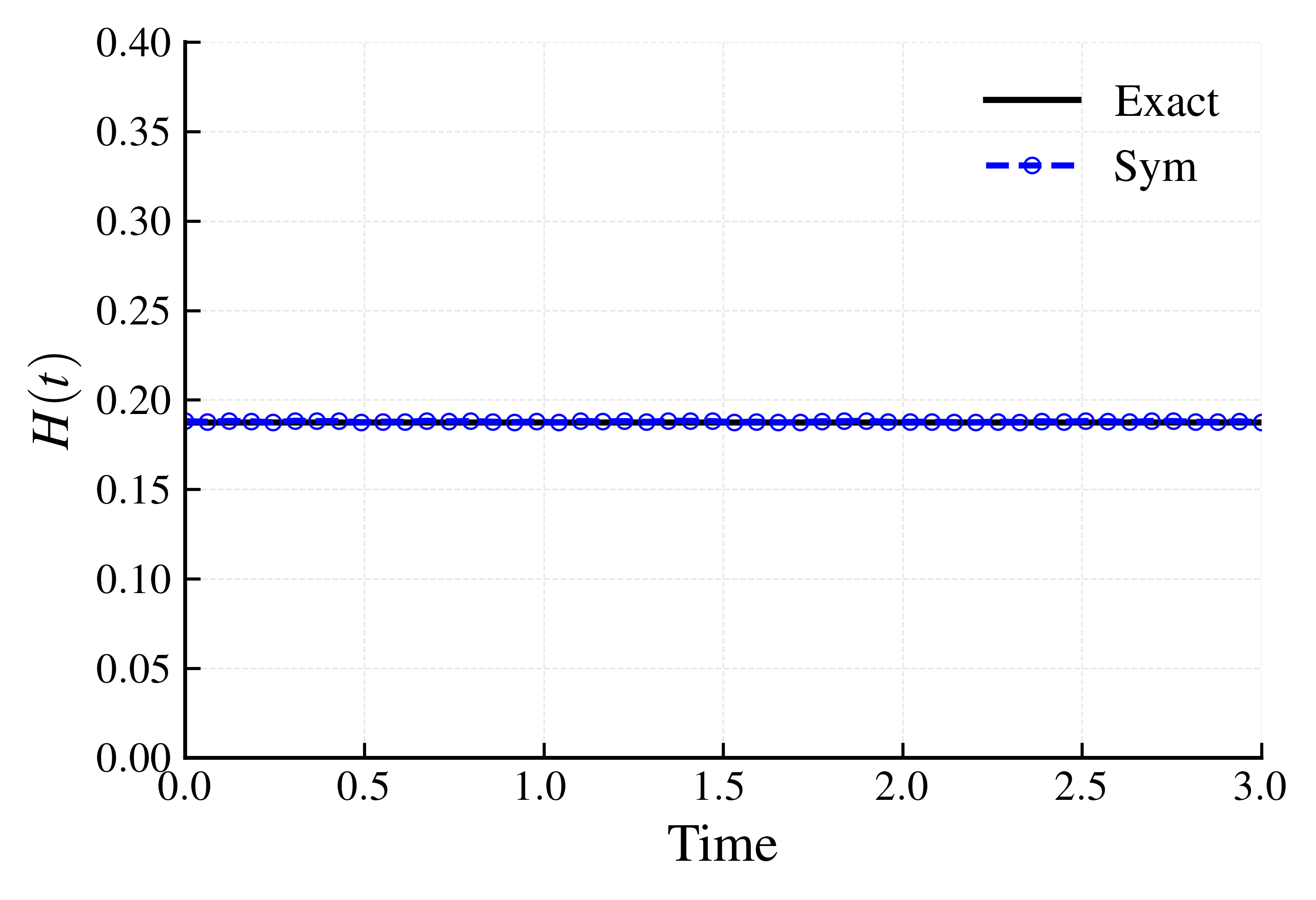}}\quad\quad
     \subfigure[Hamiltonian deviation]{\includegraphics[width=0.45\textwidth]{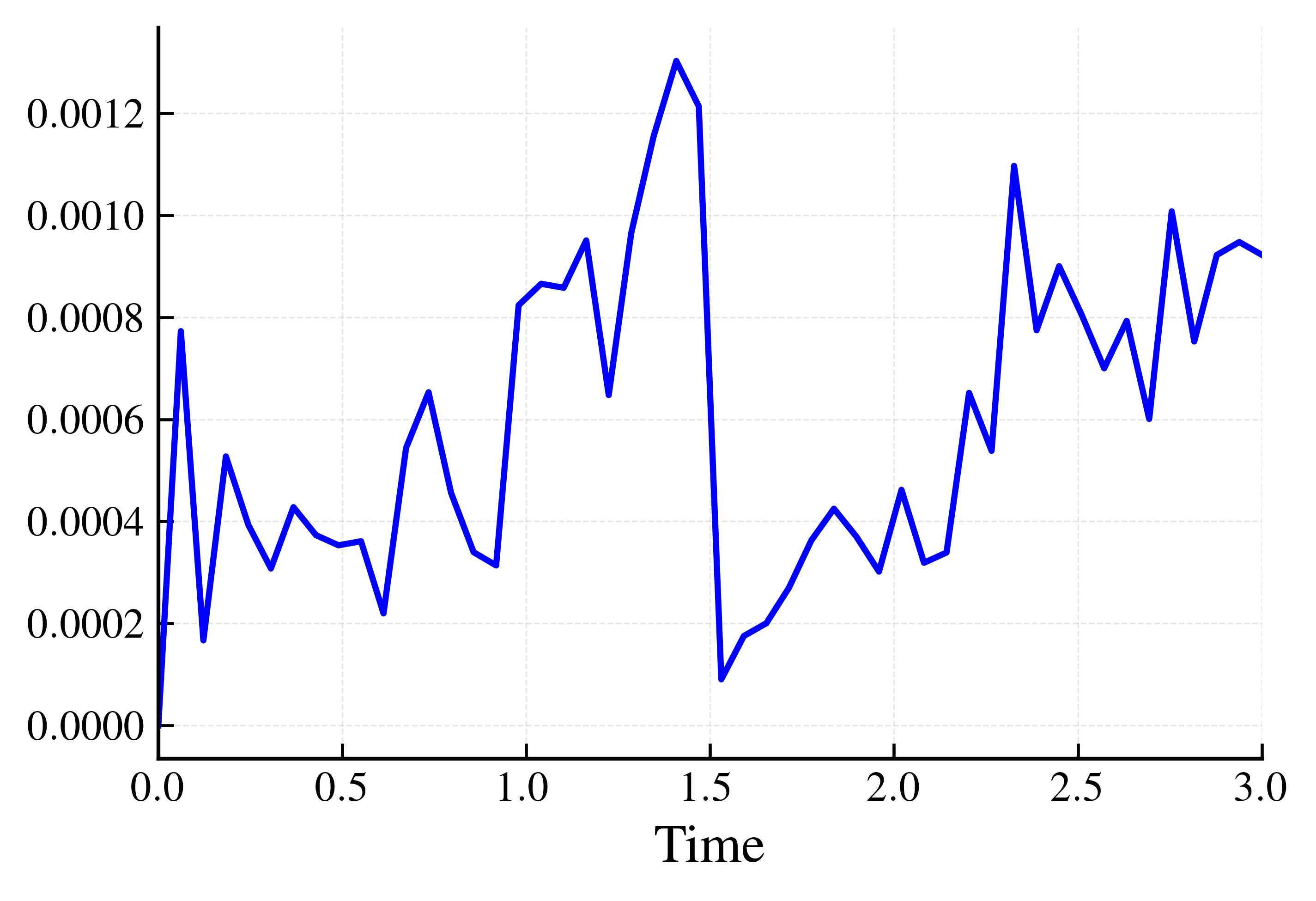}}
	\caption{Hamiltonian conservation results for a single trajectory in Example~\ref{ex:param-linear}: (a) discrete Hamiltonian in time, and (b) instantaneous Hamiltonian deviation from the initial value. }\label{fig:param_ham}
\end{figure} 

Figure~\ref{fig:param_linear} shows the representative result for the parameter choice $\kappa(\omega)=0.01736$. 
Figure~\ref{fig:param_linear}(a) shows that the latent trajectories are smooth and well-behaved, 
while
Figures~\ref{fig:param_linear}(b)--(d) show the reconstructed solutions at different time instances. The results indicate that our symplectic method successfully captures the dynamics of the full-order model accurately.
The conservation properties of the method are illustrated in Figure~\ref{fig:param_ham}, which shows the temporal evolution of the discrete Hamiltonian. The discrete Hamiltonian remains nearly constant with only minor fluctuations around the initial value, confirming the structure-preserving nature of our approach. 

\begin{table}[H]
\centering
\caption{Test reconstruction MSE in Example~\ref{ex:param-linear}, including H-OpInf baseline. For trainable neural architectures, the reported values are mean $\pm$ standard deviation over three independent runs; the cotangent-lift and H-OpInf results are deterministic.}
\label{tab:param-rec}
\small
\setlength{\tabcolsep}{3pt}

\begin{tabular}{|l|c|c|c|c|c|}
\hline
Method & Cotangent-lift & G-reflector & H\'enonNet & H\'enonNet + G-reflector & H-OpInf \\ \hline
MSE & 5.263E-5 & 
\makecell{3.928E-5 \\ $\pm$ \scriptsize{2.715E-6}} & 
\makecell{6.441E-7 \\ $\pm$ \scriptsize{6.067E-8}} & 
\makecell{6.267E-7 \\ $\pm$ \scriptsize{4.434E-8}} & 
6.832E-6 \\ \hline
\end{tabular}
\end{table}

A comparison of reconstruction accuracy across five approaches is presented: the cotangent-lift method, a 100-layer G-reflector implementation, a standalone H\'enonNet architecture, a composite H\'enonNet + G-reflector architecture, and the H-OpInf baseline \cite{gruber2025variationally}. For the H-OpInf baseline, multiple H-OpInf-style basis families were evaluated, with the best reconstruction performance achieved using the centered POD basis. Table~\ref{tab:param-rec} reports the MSEs of reconstruction over all test trajectories. The results demonstrate that the H\'enonNet-based methods achieve the highest fidelity, with the standalone and composite architectures performing nearly identically. The H-OpInf baseline, while preserving Hamiltonian structure, exhibits a larger reconstruction MSE than all symplectic ROM configurations considered here. These results indicate that, for this parametric example, the nonlinear symplectic latent-space representation already captures the solution manifold effectively, with only marginal additional benefit from the G-reflector corrections.

\end{example}

\begin{example}\label{ex:schro}
    As the final example of the numerical section, we consider the one-dimensional, nonlinear, parametric Schr{\"o}dinger equation,
\begin{equation}\label{eq:schro}
    \left\{\begin{array}{l}
iu_{t}(x,t,\epsilon) = -u_{xx}(x,t,\epsilon) - \epsilon |u(x,t,\epsilon)|^2 u(x,t,\epsilon),\\
u(x,0)= u_0(x),\\
\end{array}\right.
\end{equation}
where $u$ is a complex valued wave function, $|\cdot|$ is the modulus operator and $\epsilon \in [0.9,1.1]$. We consider
 periodic boundary conditions; i.e., $x$ belongs to a one-dimensional torus of length $L$. The initial condition is given by
 \[
 u_0(x) = \frac{\sqrt{2}}{\cosh (x-x_0)}\exp \left(i\frac{c(x-x_0)}{x}\right), 
 \]
 where $c$ denotes the wave speed, and $x_0$ is the initial position offset. When $\epsilon=1$, $|u(x,t)|$ forms a solitary wave propagating with speed $c$; for other values, the solution consists of interacting solitary waves moving bidirectionally~\cite{sulem2007nonlinear}.

Introducing the real and imaginary variables via $u = p+ i q$, we obtain the canonical form
\begin{equation}\label{eq:schro-can}
    \left\{\begin{array}{l}
q_{t}=p_{xx}+\epsilon (q^2+p^2)p \\
p_{t}=-q_{xx} -\epsilon(q^2+p^2)q
\end{array}\right.
\end{equation}
with the continuous Hamiltonian given by
\begin{equation}
    H_{\rm cont}(q,p) = \int_0^L(q_x^2+p_x^2)+\epsilon/2(q^2+p^2)^2dx
\end{equation}

Numerically, we discretize the spatial domain into $N$ equidistant points $x_i=i\Delta x$, $\Delta x=L/N$. Defining discrete variables $q_i=q(t,x_i,\epsilon)$, $p_i=p(t,x_i,\epsilon)$, and state vector $\mathbf{z}=(q_1,\dots,q_N,p_1,\dots,p_N)$, the system~\eqref{eq:schro-can} is discretized using second-order central finite difference,
 \begin{equation} 
     \frac{d}{d t} \mathbf{z}=\bbJ_{2 N} \nabla_{\bf z} H({\bf z})=\bbJ_{2 N} L\mathbf{z} + \bbJ_{2 N}\mathbf{g}(\bf z),
 \end{equation}
 where
 \begin{equation}
     L=\left(\begin{array}{cc}
D_{x x} & 0_{N} \\
0_{N} &  D_{x x}
\end{array}\right),
 \end{equation}
and $\mathbf{g}$ is a vector valued nonlinear function defined as
\begin{equation}
   \mathbf{g}(\bf z) =\left(\begin{array}{c}
(q_1^2+p_1^2)q_1  \\
\vdots  \\
(q_N^2+p_N^2)q_N \\
(q_1^2+p_1^2)p_1  \\
\vdots  \\
(q_N^2+p_N^2)p_N 
\end{array}\right).
\end{equation}
The discrete Hamiltonian is
\begin{equation}\label{eq:schro-ham-discre}
    H(\mathbf{z}) = \Delta x\sum_{i=1}^{N}\left(  \frac{\left(q_{i-1}q_i-q_{i}^2\right)^{2}}{\Delta x^{2}}+\frac{\left(p_{i-1}p_i-p_{i}^2\right)^{2}}{\Delta x^{2}}+\frac{\epsilon}{4}(p_i^2+q_i^2) \right).
\end{equation}



\begin{table}[H]
\centering
\caption{Numerical discretization parameters used in Example~\ref{ex:schro}.}
\label{tab:ex5-3-config}
\small
\begin{tabular}{|l|c|}
\hline
Domain length              & $L = 2\pi/l$ \\ \hline
Domain scaling factor      & $l = 0.11$ \\ \hline
Grid points                & $N = 256$ \\ \hline
Time discretization size   & $\Delta t = 0.2$ \\ \hline
Wave speed                 & $c = 1$ \\ \hline
\end{tabular}
\end{table}

\begin{table}[H]
\centering
\caption{Neural network architectures used in Example~\ref{ex:schro}.}
\label{tab:ex5-3-network}
\small
\begin{tabular}{|l|c|}
\hline
Autoencoder H\'enonNet      & $[512]\times 2$ \\ \hline
Autoencoder G-reflector     & $10$ \\ \hline
Flow mapping H\'enonNet     & $[128]\times 2$ \\ \hline
\end{tabular}
\end{table}

Because the semi-discrete Hamiltonian is nonlinear, the St{\"o}rmer–Verlet time integrator must be applied in its implicit form.  At each time step the resulting nonlinear system is resolved by a Newton–Raphson iteration, which is terminated once the residual infinity-norm falls below an absolute tolerance of $10^{-12}$. 
Snapshot data are generated by uniformly sampling 500 parameter values $\{\epsilon_i\}_{i=1}^{500}\subset[0.9,1.1]$, integrated numerically over $[0,20]$. In addition, we randomly sample 100 parameter values from the same parameter range as an independent test set for reconstruction evaluation. The reduced model employs a 12-dimensional latent space. Configuration and network architecture details are provided in Tables~\ref{tab:ex5-3-config} and~\ref{tab:ex5-3-network}, with hyperparameters $(\lambda_1,\lambda_2)$ in loss~\eqref{eq:loss_total} set to be $(1,0.05)$. 

\begin{figure}[!htbp]
 \centering
     \subfigure[Latent space dynamics]{\includegraphics[width=0.45\textwidth]{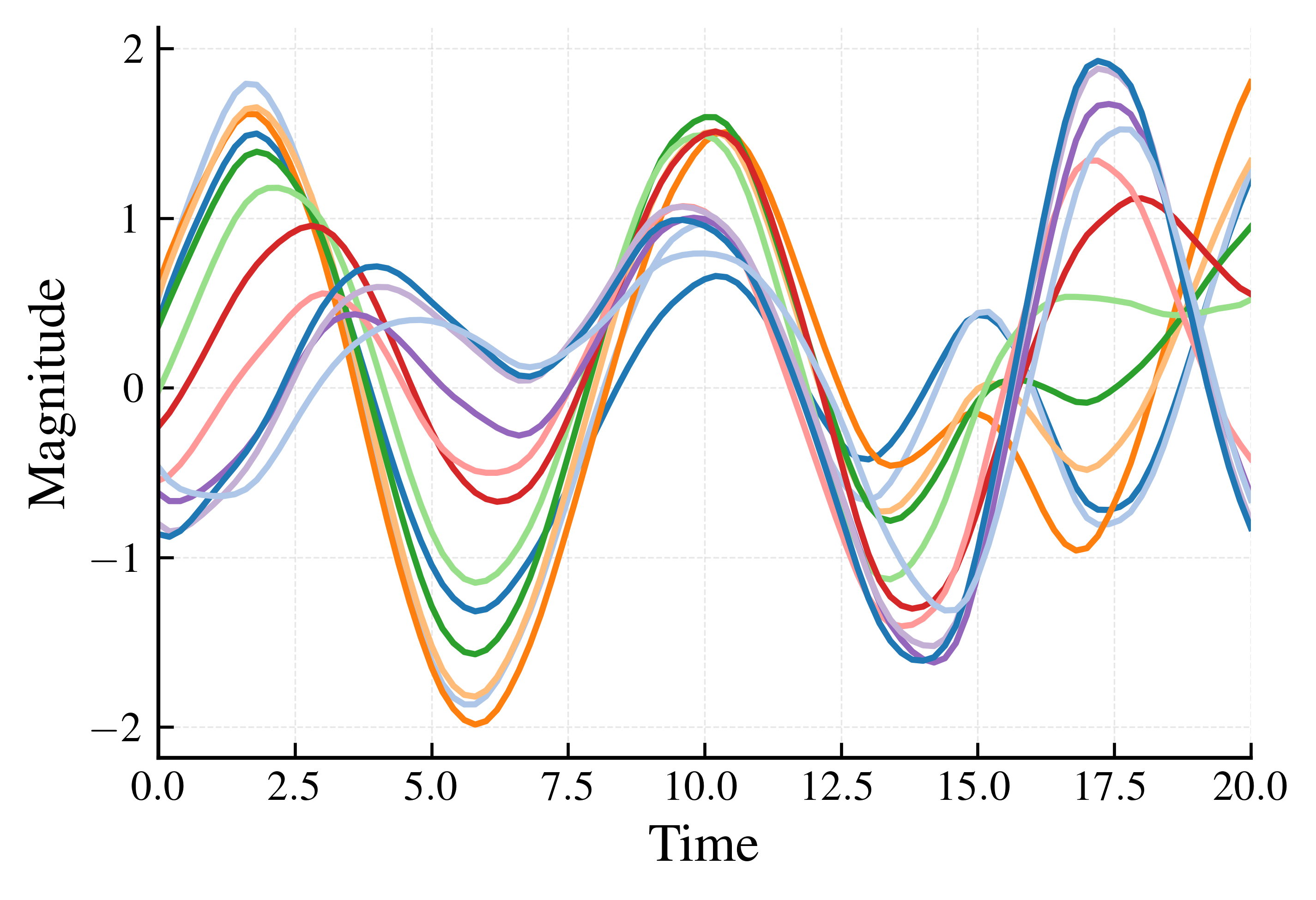}}\quad
     \subfigure[$t=0$]{\includegraphics[width=0.45\textwidth]{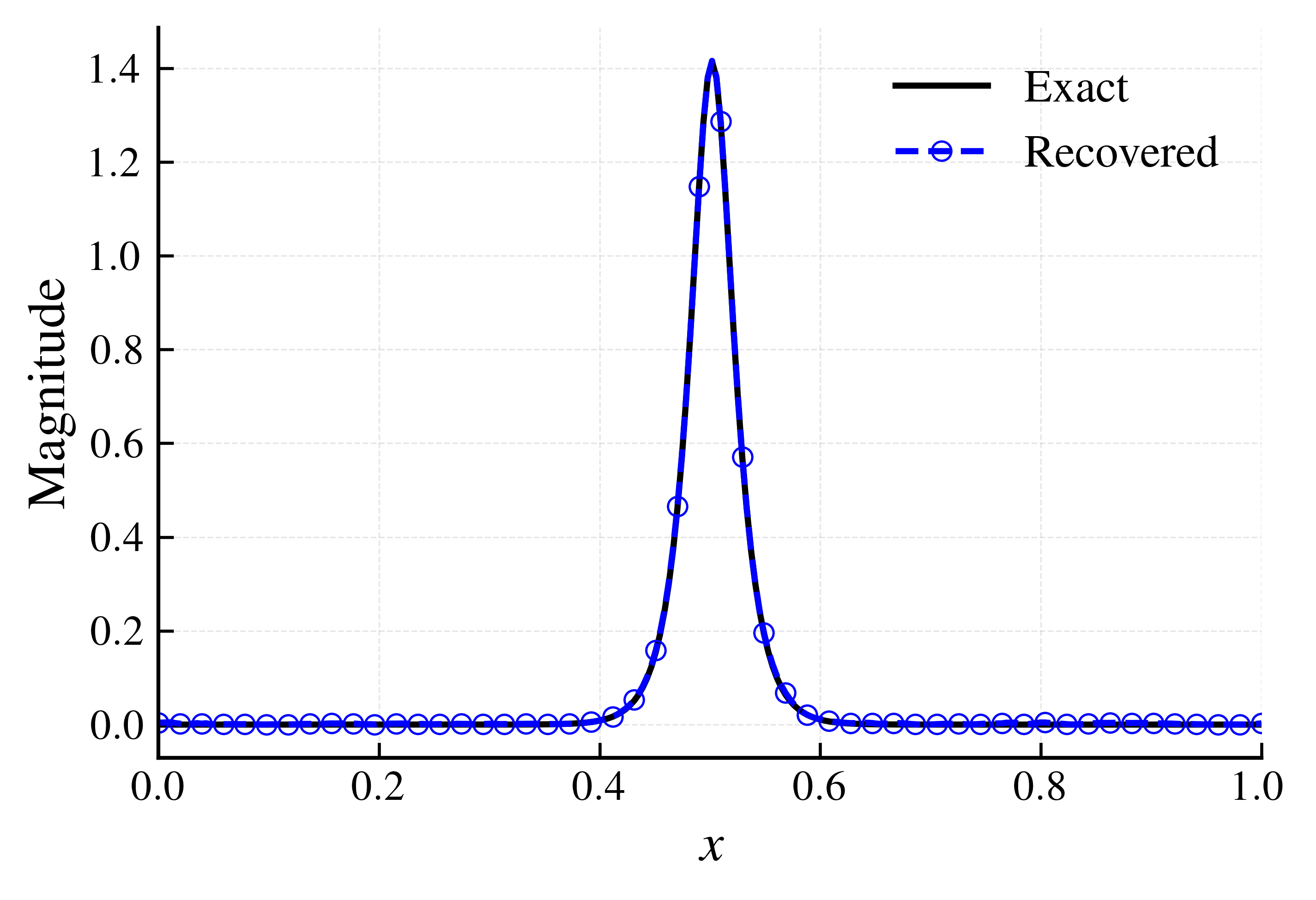}}\\
     \subfigure[$t=10$]{\includegraphics[width=0.45\textwidth]{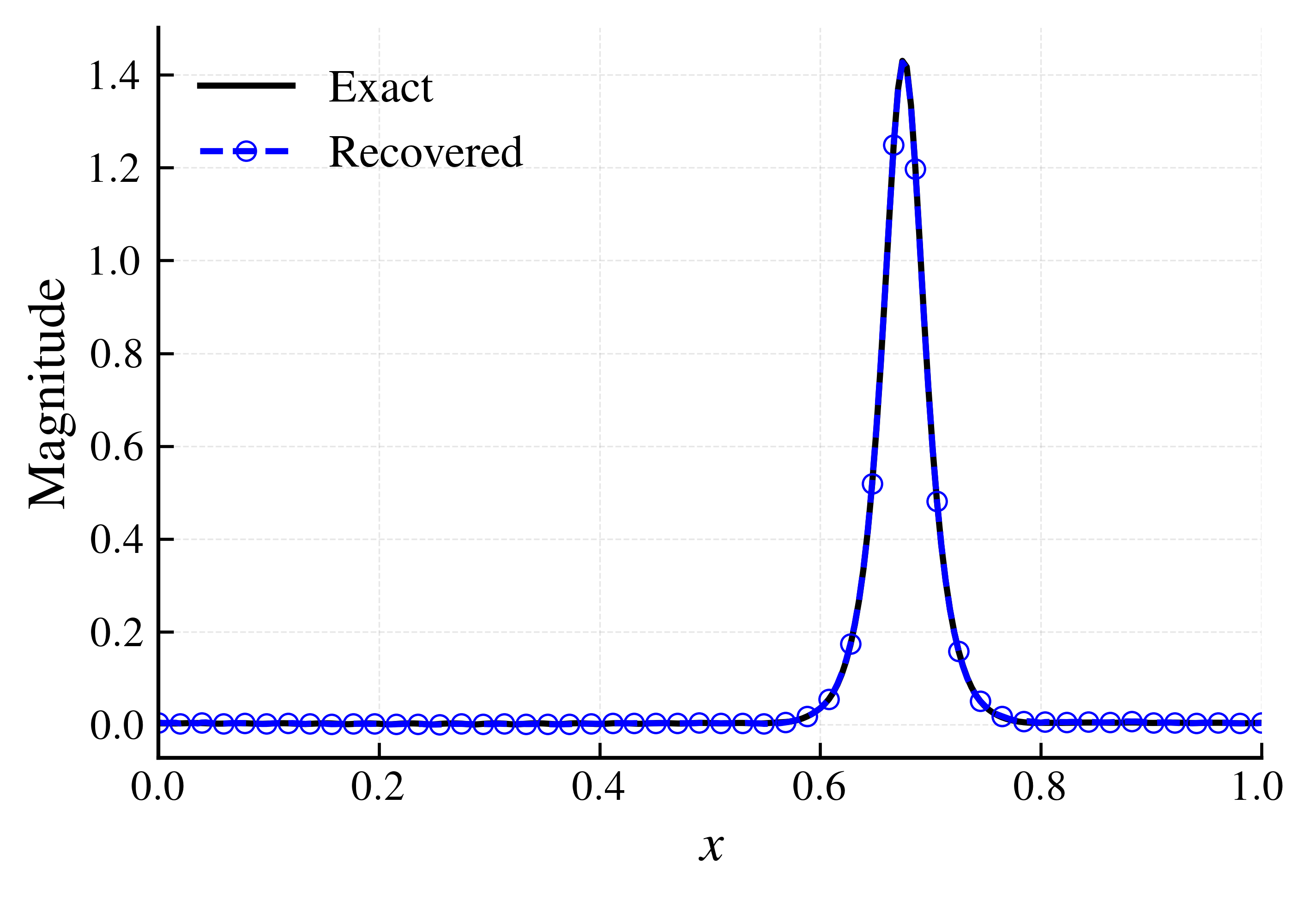}}\quad
     \subfigure[$t=20$]{\includegraphics[width=0.45\textwidth]{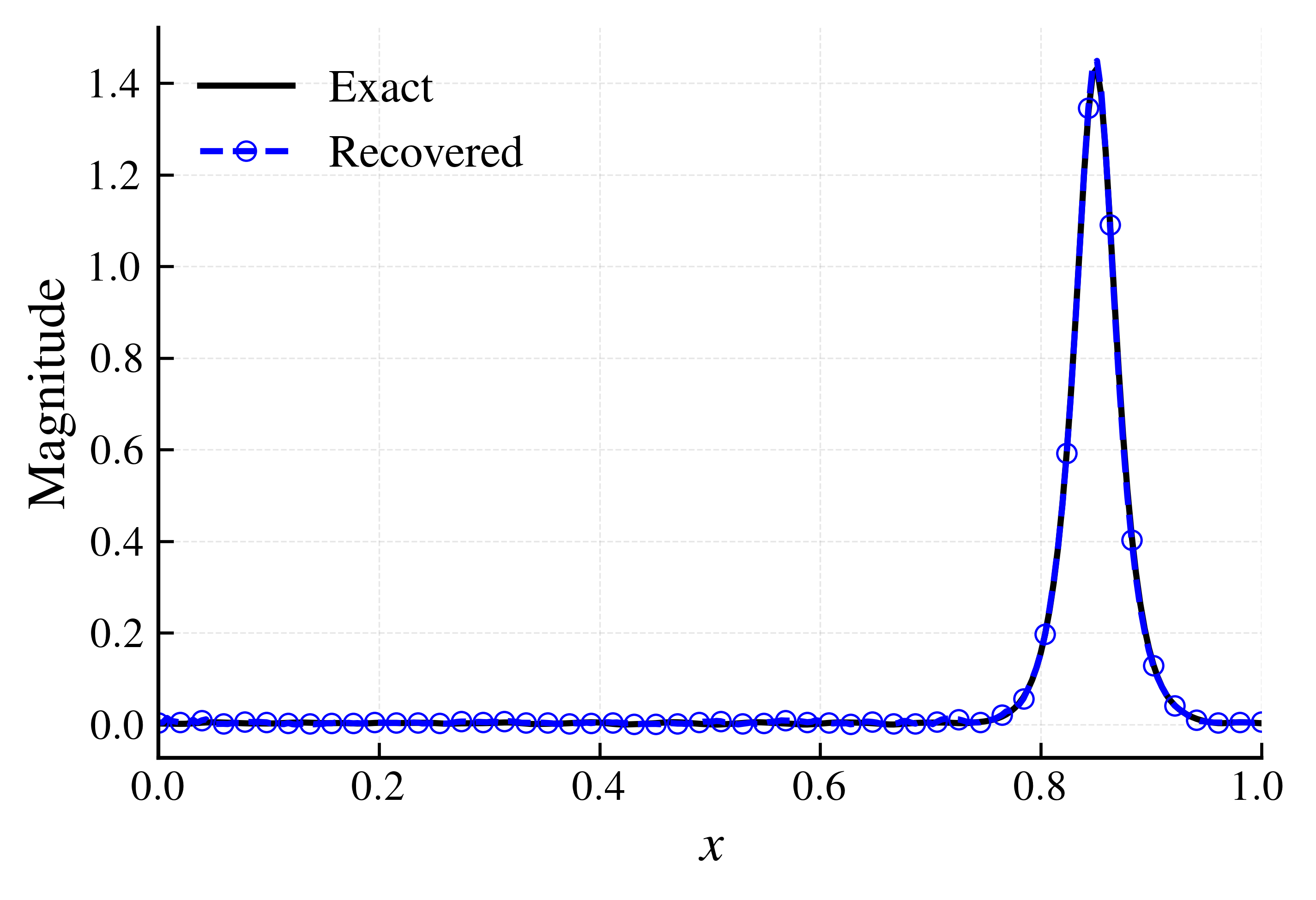}}
	\caption{ROM prediction in Example~\ref{ex:schro} with $\epsilon=1.007858$: (a) latent space dynamics, and reconstructed solutions at times (b) $t=0$, (c) $t=10$, and (d) $t=20$.}\label{fig:schro}
\end{figure}

Figure \ref{fig:schro}(a) illustrates the latent space trajectories for a representative case with parameter value $\epsilon = 1.007858$. The corresponding reconstructed solutions at times $t = 0$, $t = 10$, and $t = 20$ are shown in Figures \ref{fig:schro}(b)-(d). Despite the system's inherent nonlinearity, our method generates a stable reduced-order model that accurately approximates the full system dynamics over the entire integration period.

\begin{figure}[!htbp]
 \centering
     \subfigure[Hamiltonian]{\includegraphics[width=0.45\textwidth]{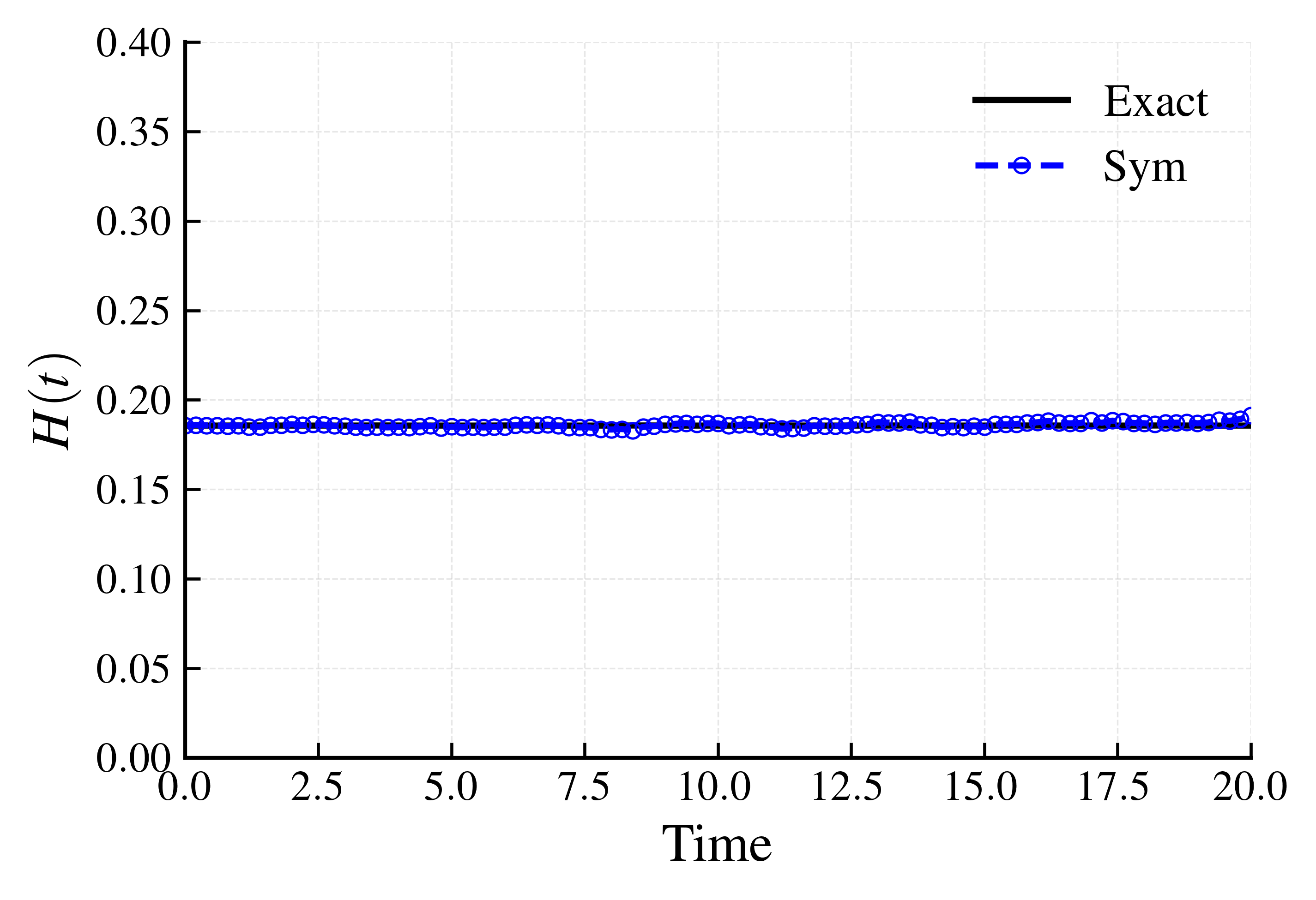}}\quad\quad
     \subfigure[Hamiltonian deviation]{\includegraphics[width=0.45\textwidth]{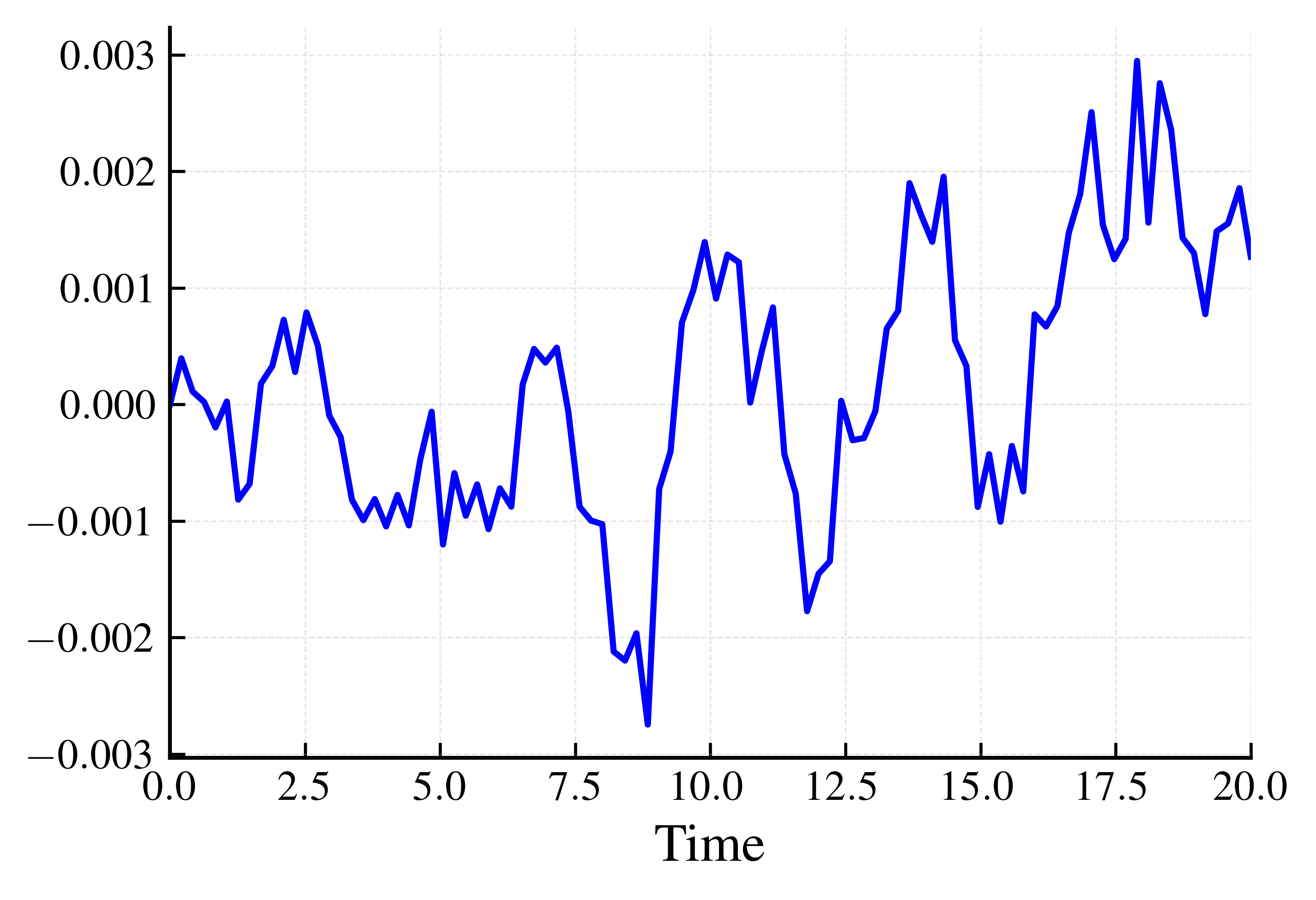}}
	\caption{Hamiltonian conservation results for a single trajectory in Example~\ref{ex:schro}: (a) discrete Hamiltonian in time, and (b) instantaneous Hamiltonian deviation from the initial value. }\label{fig:schro-ham}
\end{figure}

The conservation properties of our symplectic approach are demonstrated in Figure~\ref{fig:schro-ham}, which displays both the temporal evolution of the discrete Hamiltonian and the instantaneous Hamiltonian deviation from its initial value. The Hamiltonian remains nearly constant throughout the simulation, with only very small deviations from its initial value, confirming that our method effectively preserves the underlying Hamiltonian structure even in this nonlinear setting.

\begin{table}[H]
\centering
\caption{Test reconstruction MSE in Example~\ref{ex:schro}, including H-OpInf baseline. For trainable neural architectures, the reported values are mean $\pm$ standard deviation over three independent runs; the cotangent-lift and H-OpInf results are deterministic.}
\label{tab:schro-rec}
\small
\setlength{\tabcolsep}{3pt}

\begin{tabular}{|l|c|c|c|c|c|}
\hline
Method & Cotangent-lift & G-reflector & H\'enonNet & H\'enonNet + G-reflector & H-OpInf \\ \hline
MSE & 6.175E-4 & 
\makecell{2.029E-4 \\ $\pm$ \scriptsize{1.435E-5}} & 
\makecell{6.378E-7 \\ $\pm$ \scriptsize{6.110E-8}} & 
\makecell{6.128E-7 \\ $\pm$ \scriptsize{4.130E-8}} & 
2.331E-4 \\ \hline
\end{tabular}
\end{table}

Table~\ref{tab:schro-rec} presents the test reconstruction MSEs for five approaches: the cotangent-lift method, a 100-layer G-reflector implementation, a standalone H\'enonNet architecture, a composite H\'enonNet + G-reflector architecture, and the H-OpInf baseline. As with the previous examples, the H-OpInf baseline was evaluated over multiple basis families, and the best performance was obtained with the centered POD basis. The H\'enonNet-based methods achieve substantially lower reconstruction errors than the linear embeddings and the H-OpInf baseline, with the standalone and composite H\'enonNet architectures again performing nearly identically. The larger performance gap observed in this strongly nonlinear example suggests that the proposed nonlinear symplectic ROM is particularly effective when the underlying solution manifold exhibits pronounced nonlinear phase-space coupling. This observation is consistent with the ability of H\'enonNet-based mappings to model such nonlinear coupling more effectively in the setting considered here.


\end{example}

\section{Conclusion}\label{sec:conclusion}
In this paper, we have introduced a novel data-driven symplectic reduced-order modeling (ROM) framework for Hamiltonian systems that unifies dimensionality reduction and dynamics learning within a single, end-to-end architecture. The method is built on H\'enon neural networks (H\'enonNets) as the primary nonlinear blocks, optionally augmented with linear 
G-reflector layers for efficient linear symplectic corrections. Through three Hamiltonian examples—the linear wave equation, a parametric linear wave equation, and the nonlinear Schr{\"o}dinger equation—we validated the performance and robustness of our framework with increasing levels of complexity, from linear and parametric scenarios to highly nonlinear systems. 
Numerical results demonstrate that our method accurately reconstructs the original high-dimensional dynamics, maintains long-term predictive capability even beyond the training horizon, and preserves the Hamiltonian structure, as evidenced by minimal fluctuation in the discrete Hamiltonian.
Overall, the proposed data-driven symplectic ROM approach shows promising results for model reduction of Hamiltonian systems. 

One limitation of the proposed work is its data-driven nature, which constrains the incorporation of additional physical information. An intrusive ROM approach, built on the proposed structure-preserving compression, would therefore be an interesting direction to explore.
Another future direction is to extend this framework to broader classes of nonlinear PDEs and higher-dimensional dynamical systems, such as particle dynamics and their collective effects in accelerators \cite{huang2024symplectic}.

\bibliographystyle{siamplain}
\bibliography{references}
\end{document}

%% file: ex_shared.tex
\usepackage{lipsum}
\usepackage{amsfonts}
\usepackage{makecell}   
\usepackage{amsmath}
\usepackage{graphicx}
\usepackage{epstopdf}
\usepackage{algorithmic}
\newtheorem{remark}{Remark}
\usepackage{multirow}
\usepackage{subfigure}
\newtheorem{example}{Example}[section]

\usepackage{soul}

\newcommand{\xv}{\mathbf{ x}}
\newcommand{\yv}{\mathbf{ y}}

\newcommand{\grad}{\nabla}
\newcommand{\bbJ}{\mathbb{J}}

\newcommand{\fdec}{f_{\rm dec}}
\newcommand{\fenc}{f_{\rm enc}}
\newcommand{\fflow}{f_{\rm flow}}
\ifpdf
  \DeclareGraphicsExtensions{.eps,.pdf,.png,.jpg}
\else
  \DeclareGraphicsExtensions{.eps}
\fi

\newsiamremark{hypothesis}{Hypothesis}
\crefname{hypothesis}{Hypothesis}{Hypotheses}
\newsiamthm{claim}{Claim}
\newsiamremark{fact}{Fact}
\crefname{fact}{Fact}{Facts}

\headers{Reduced-order modeling of Hamiltonian dynamics}{Y. Chen, W. Guo, Q. Tang, and X. Zhong}

\title{Reduced-order modeling of Hamiltonian dynamics based on symplectic neural networks\thanks{Submitted to the editors August 2025.
\funding{This work was partially supported by the U.S. Department of Energy Advanced Scientific Computing Research (ASCR) program's Mathematical Multifaceted Integrated Capability Center (MMICC), 
the ASCR SciML program under DE-FOA-2493 “Data-intensive scientific machine learning”, 
{Department of
Energy Office of Science, Early Career Research Program under Award Number DE-SC0026277,}
and the Air Force Office of Scientific Research FA9550-18-1-0257.}}}
\author{
Yongsheng Chen\thanks{School of Mathematical Sciences, Zhejiang University, Hangzhou 310027, China. (\email{22035024@zju.edu.cn, zhongxh@zju.edu.cn}).}
\and Wei Guo\thanks{Department of Mathematics and Statistics, Texas Tech University, Lubbock, TX 70409, USA 
  (\email{weimath.guo@ttu.edu}).}
\and Qi Tang\thanks{School of Computational Science and Engineering, Georgia Institute of Technology, Atlanta, GA 30332, USA.}
\and Xinghui Zhong\footnotemark[2]
}

\usepackage{amsopn}
